\documentclass[11pt,twoside,leqno]{aomamlt2e} 
  \pageno{879}
\received{February 21, 2003}
 \newcounter{notes}

\numberwithin{equation}{section}

 \newcommand{\ee}{{\hskip1pt\rm \'{\hskip-6.5pt \it e}}}
\newcommand\hfq{\hfill\qed}

 \newtheorem{theorem}{Theorem}[section]
\newtheorem{lem}[theorem]{Lemma}
\newtheorem{cor}[theorem]{Corollary}
\theorembodyfont{\upshape}
\newtheorem{defi}[theorem]{{\it Definition}}
\newtheorem{rem}[theorem]{{\it Remark}}
 
   \catcode`\@=11
\def\relaxnext@{\let\next\relax}
\font\tenmsx=msam10       
\font\sevenmsx=msam7      
\font\fivemsx=msam5       
\font\tenmsy=msbm10       
\font\sevenmsy=msbm7      
\font\fivemsy=msbm5       
\newfam\msxfam
\newfam\msyfam
\textfont\msxfam=\tenmsx  \scriptfont\msxfam=\sevenmsx
  \scriptscriptfont\msxfam=\fivemsx
\textfont\msyfam=\tenmsy  \scriptfont\msyfam=\sevenmsy
  \scriptscriptfont\msyfam=\fivemsy

\def\hexnumber@#1{\ifcase#1 0\or1\or2\or3\or4\or5\or6\or7\or8\or9\or
	A\or B\or C\or D\or E\or F\fi }

\font\teneuf=eufm10
\font\seveneuf=eufm7
\font\fiveeuf=eufm5
\newfam\euffam
\textfont\euffam=\teneuf
\scriptfont\euffam=\seveneuf
\scriptscriptfont\euffam=\fiveeuf
\def\frak{\relaxnext@\ifmmode\let\next\frak@\else
 \def\next{\errmessage{Use \string\frak\space only in math
mode}}\fi\next}
\def\goth{\relaxnext@\ifmmode\let\next\frak@\else 
\def\next{\errmessage{Use \string\goth\space only in math
mode}}\fi\next} 
\def\frak@#1{{\frak@@{#1}}}
\def\frak@@#1{\fam\euffam#1}

\edef\msx@{\hexnumber@\msxfam}
\edef\msy@{\hexnumber@\msyfam}

\mathchardef\boxdot="2\msx@00
\mathchardef\boxplus="2\msx@01
\mathchardef\boxtimes="2\msx@02
\mathchardef\square="0\msx@03
\mathchardef\blacksquare="0\msx@04
\mathchardef\centerdot="2\msx@05
\mathchardef\lozenge="0\msx@06
\mathchardef\blacklozenge="0\msx@07
\mathchardef\circlearrowright="3\msx@08
\mathchardef\circlearrowleft="3\msx@09
\mathchardef\rightleftharpoons="3\msx@0A
\mathchardef\leftrightharpoons="3\msx@0B
\mathchardef\boxminus="2\msx@0C
\mathchardef\Vdash="3\msx@0D
\mathchardef\Vvdash="3\msx@0E
\mathchardef\vDash="3\msx@0F
\mathchardef\twoheadrightarrow="3\msx@10
\mathchardef\twoheadleftarrow="3\msx@11
\mathchardef\leftleftarrows="3\msx@12
\mathchardef\rightrightarrows="3\msx@13
\mathchardef\upuparrows="3\msx@14
\mathchardef\downdownarrows="3\msx@15
\mathchardef\upharpoonright="3\msx@16

\mathchardef\downharpoonright="3\msx@17
\mathchardef\upharpoonleft="3\msx@18
\mathchardef\downharpoonleft="3\msx@19
\mathchardef\rightarrowtail="3\msx@1A
\mathchardef\leftarrowtail="3\msx@1B
\mathchardef\leftrightarrows="3\msx@1C
\mathchardef\rightleftarrows="3\msx@1D
\mathchardef\Lsh="3\msx@1E
\mathchardef\Rsh="3\msx@1F
\mathchardef\rightsquigarrow="3\msx@20
\mathchardef\leftrightsquigarrow="3\msx@21
\mathchardef\looparrowleft="3\msx@22
\mathchardef\looparrowright="3\msx@23
\mathchardef\circeq="3\msx@24
\mathchardef\succsim="3\msx@25
\mathchardef\gtrsim="3\msx@26
\mathchardef\gtrapprox="3\msx@27
\mathchardef\multimap="3\msx@28
\mathchardef\therefore="3\msx@29
\mathchardef\because="3\msx@2A
\mathchardef\doteqdot="3\msx@2B

\mathchardef\triangleq="3\msx@2C
\mathchardef\precsim="3\msx@2D
\mathchardef\lesssim="3\msx@2E
\mathchardef\lessapprox="3\msx@2F
\mathchardef\eqslantless="3\msx@30
\mathchardef\eqslantgtr="3\msx@31
\mathchardef\curlyeqprec="3\msx@32
\mathchardef\curlyeqsucc="3\msx@33
\mathchardef\preccurlyeq="3\msx@34
\mathchardef\leqq="3\msx@35
\mathchardef\leqslant="3\msx@36
\mathchardef\lessgtr="3\msx@37
\mathchardef\backprime="0\msx@38
\mathchardef\risingdotseq="3\msx@3A
\mathchardef\fallingdotseq="3\msx@3B
\mathchardef\succcurlyeq="3\msx@3C
\mathchardef\geqq="3\msx@3D
\mathchardef\geqslant="3\msx@3E
\mathchardef\gtrless="3\msx@3F
\mathchardef\sqsubset="3\msx@40
\mathchardef\sqsupset="3\msx@41
\mathchardef\vartriangleright="3\msx@42
\mathchardef\vartriangleleft="3\msx@43
\mathchardef\trianglerighteq="3\msx@44
\mathchardef\trianglelefteq="3\msx@45
\mathchardef\bigstar="0\msx@46
\mathchardef\between="3\msx@47
\mathchardef\blacktriangledown="0\msx@48
\mathchardef\blacktriangleright="3\msx@49
\mathchardef\blacktriangleleft="3\msx@4A
\mathchardef\vartriangle="0\msx@4D
\mathchardef\blacktriangle="0\msx@4E
\mathchardef\triangledown="0\msx@4F
\mathchardef\lesseqgtr="3\msx@51
\mathchardef\gtreqless="3\msx@52
\mathchardef\lesseqqgtr="3\msx@53
\mathchardef\gtreqqless="3\msx@54
\mathchardef\Rrightarrow="3\msx@56
\mathchardef\Lleftarrow="3\msx@57
\mathchardef\veebar="2\msx@59
\mathchardef\barwedge="2\msx@5A
\mathchardef\doublebarwedge="2\msx@5B
\mathchardef\angle="0\msx@5C
\mathchardef\measuredangle="0\msx@5D
\mathchardef\sphericalangle="0\msx@5E
\mathchardef\varpropto="3\msx@5F
\mathchardef\smallsmile="3\msx@60
\mathchardef\smallfrown="3\msx@61
\mathchardef\Subset="3\msx@62
\mathchardef\Supset="3\msx@63
\mathchardef\Cup="2\msx@64

\mathchardef\Cap="2\msx@65

\mathchardef\curlywedge="2\msx@66
\mathchardef\curlyvee="2\msx@67
\mathchardef\leftthreetimes="2\msx@68
\mathchardef\rightthreetimes="2\msx@69
\mathchardef\subseteqq="3\msx@6A
\mathchardef\supseteqq="3\msx@6B
\mathchardef\bumpeq="3\msx@6C
\mathchardef\Bumpeq="3\msx@6D
\mathchardef\lll="3\msx@6E

\mathchardef\ggg="3\msx@6F

\mathchardef\circledS="0\msx@73
\mathchardef\pitchfork="3\msx@74
\mathchardef\dotplus="2\msx@75
\mathchardef\backsim="3\msx@76
\mathchardef\backsimeq="3\msx@77
\mathchardef\complement="0\msx@7B
\mathchardef\intercal="2\msx@7C
\mathchardef\circledcirc="2\msx@7D
\mathchardef\circledast="2\msx@7E
\mathchardef\circleddash="2\msx@7F
\def\ulcorner{\delimiter"4\msx@70\msx@70 }
\def\urcorner{\delimiter"5\msx@71\msx@71 }
\def\llcorner{\delimiter"4\msx@78\msx@78 }
\def\lrcorner{\delimiter"5\msx@79\msx@79 }
\def\yen{\mathhexbox\msx@55 }
\def\checkmark{\mathhexbox\msx@58 }
\def\circledR{\mathhexbox\msx@72 }
\def\maltese{\mathhexbox\msx@7A }
\mathchardef\lvertneqq="3\msy@00
\mathchardef\gvertneqq="3\msy@01
\mathchardef\nleq="3\msy@02
\mathchardef\ngeq="3\msy@03
\mathchardef\nless="3\msy@04
\mathchardef\ngtr="3\msy@05
\mathchardef\nprec="3\msy@06
\mathchardef\nsucc="3\msy@07
\mathchardef\lneqq="3\msy@08
\mathchardef\gneqq="3\msy@09
\mathchardef\nleqslant="3\msy@0A
\mathchardef\ngeqslant="3\msy@0B
\mathchardef\lneq="3\msy@0C
\mathchardef\gneq="3\msy@0D
\mathchardef\npreceq="3\msy@0E
\mathchardef\nsucceq="3\msy@0F
\mathchardef\precnsim="3\msy@10
\mathchardef\succnsim="3\msy@11
\mathchardef\lnsim="3\msy@12
\mathchardef\gnsim="3\msy@13
\mathchardef\nleqq="3\msy@14
\mathchardef\ngeqq="3\msy@15
\mathchardef\precneqq="3\msy@16
\mathchardef\succneqq="3\msy@17
\mathchardef\precnapprox="3\msy@18
\mathchardef\succnapprox="3\msy@19
\mathchardef\lnapprox="3\msy@1A
\mathchardef\gnapprox="3\msy@1B
\mathchardef\nsim="3\msy@1C
\mathchardef\ncong="3\msy@1D

\mathchardef\varsubsetneq="3\msy@20
\mathchardef\varsupsetneq="3\msy@21
\mathchardef\nsubseteqq="3\msy@22
\mathchardef\nsupseteqq="3\msy@23
\mathchardef\subsetneqq="3\msy@24
\mathchardef\supsetneqq="3\msy@25
\mathchardef\varsubsetneqq="3\msy@26
\mathchardef\varsupsetneqq="3\msy@27
\mathchardef\subsetneq="3\msy@28
\mathchardef\supsetneq="3\msy@29
\mathchardef\nsubseteq="3\msy@2A
\mathchardef\nsupseteq="3\msy@2B
\mathchardef\nparallel="3\msy@2C
\mathchardef\nmid="3\msy@2D
\mathchardef\nshortmid="3\msy@2E
\mathchardef\nshortparallel="3\msy@2F
\mathchardef\nvdash="3\msy@30
\mathchardef\nVdash="3\msy@31
\mathchardef\nvDash="3\msy@32
\mathchardef\nVDash="3\msy@33
\mathchardef\ntrianglerighteq="3\msy@34
\mathchardef\ntrianglelefteq="3\msy@35
\mathchardef\ntriangleleft="3\msy@36
\mathchardef\ntriangleright="3\msy@37
\mathchardef\nleftarrow="3\msy@38
\mathchardef\nrightarrow="3\msy@39
\mathchardef\nLeftarrow="3\msy@3A
\mathchardef\nRightarrow="3\msy@3B
\mathchardef\nLeftrightarrow="3\msy@3C
\mathchardef\nleftrightarrow="3\msy@3D
\mathchardef\divideontimes="2\msy@3E
\mathchardef\varnothing="0\msy@3F
\mathchardef\nexists="0\msy@40
\mathchardef\mho="0\msy@66
\mathchardef\eth="0\msy@67
\mathchardef\eqsim="3\msy@68
\mathchardef\beth="0\msy@69
\mathchardef\gimel="0\msy@6A
\mathchardef\daleth="0\msy@6B
\mathchardef\lessdot="3\msy@6C
\mathchardef\gtrdot="3\msy@6D
\mathchardef\ltimes="2\msy@6E
\mathchardef\rtimes="2\msy@6F
\mathchardef\shortmid="3\msy@70
\mathchardef\shortparallel="3\msy@71
\mathchardef\smallsetminus="2\msy@72
\mathchardef\thicksim="3\msy@73
\mathchardef\thickapprox="3\msy@74
\mathchardef\approxeq="3\msy@75
\mathchardef\succapprox="3\msy@76
\mathchardef\precapprox="3\msy@77
\mathchardef\curvearrowleft="3\msy@78
\mathchardef\curvearrowright="3\msy@79
\mathchardef\digamma="0\msy@7A
\mathchardef\varkappa="0\msy@7B
\mathchardef\hslash="0\msy@7D
\mathchardef\hbar="0\msy@7E
\mathchardef\backepsilon="3\msy@7F
\def\Bbb{\ifmmode\let\next\Bbb@\else
 \def\next{\errmessage{Use \string\Bbb\space only in math mode}}\fi\next}
\def\Bbb@#1{{\Bbb@@{#1}}}
\def\Bbb@@#1{\fam\msyfam#1}

\catcode`\@=12

\def\complex{{\Bbb C}}
\def\real{{\Bbb R}}
\def\Natural{{\Bbb N}}

\def\coleq{\mathchar"303A=}
\def\maps{\nobreak\hskip.1111em\mathpunct{}\nonscript%
       \mkern-\thinmuskip{:}\hskip.3333emplus.0555em\relax}

\def\cA{{\cal A}}
\def\cC{{\cal C}}
\def\cF{{\cal F}}
\def\cH{{\cal H}}

\def\cL{{\cal L}}
\def\cM{{\cal M}}
\def\cO{{\cal O}}
\def\cP{{\cal P}}
\def\cQ{{\cal Q}}

\def\cT{{\cal T}}

\def\cW{{\cal W}}

\def\Fou{{\cF}}
\def\Hyp{{\cH}}
\def\Pal{{\cP}}

\def\Mer{{\cM}}
\def\Par{{\cP}}

\def\psW{{\cT}}

\def\Cusp{{\cA}}
\def\Lau{{\cL}}

\def\supp{\mathop{\rm supp}}
\def\reg{\mathop{\rm reg}}

\def\Res{\mathop{\rm Res}\nolimits}
\def\WT{\mathop{\rm WT}}

\def\Hom{\mathop{\rm Hom}\nolimits}
\def\End{\mathop{\rm End}}
\def\Span{\mathop{\rm Span}}

\def\pr{\mathop{\rm pr}\nolimits}
\def\inj{\mathop{\hbox{\rm i}}\nolimits}

\def\Re{\mathop{\rm Re}}

\def\reg{\mathop{\rm reg}}
\def\laur{{\mathop{\rm laur}}}

\def\AC{\mathop{\rm AC}\nolimits}

\def\PW{\mathop{\rm PW}\nolimits}
\def\Hil{H}

\def\rmT{{\rm T}}
\def\rmK{{\rm K}}

\def\iq{{\rm q}}

\def\bfL{{\bf L}}

\def\Parmin{\Parsigma^{{\rm min}}}

\def\after{\,{\scriptstyle\circ}\,}
\def\inp#1#2{\langle#1,#2\rangle}
\def\frac#1#2{{{#1}\over{#2}}}
\def\ay{{\fa}}

\def\Cci{C_c^\infty}
\def\Cinf{C^\infty}

 \newdimen\prestarheight\prestarheight=4pt
 \newdimen\scriptprestarheight\scriptprestarheight=2.8pt
 \newdimen\scriptscriptprestarheight\scriptscriptprestarheight=2pt
\def\prestar#1#2{\mathchoice
 {\raise#1\prestarheight\hbox{$\scriptstyle *$}\kern-.#2em}
 {\raise#1\prestarheight\hbox{$\scriptstyle *$}\kern-.#2em}
 {\raise#1\scriptprestarheight\hbox{$\scriptscriptstyle *$}\kern-.#2em}
 {\raise#1\scriptscriptprestarheight\hbox{$\scriptscriptstyle *$}\kern-.#2em}}
\def\stP{{\prestar11P}}
\def\stt{\prestar{.75}0t}

\def\Vtau{V_\tau}
\def\col{\colon}
\def\CciXt{\Cci(X\col\tau)}

\def\faqdF{\fa_{F\iq}^*}
\def\faqcdF{\fa_{F\iq\siC}^*}
\def\faqdFp{\fa^{*\perp}_{F\iq}}
\def\nE{E^\circ}
\def\nC{C^\circ}

\def\Aqp{A^+_\iq}
\def\Aq{A_\iq}
\def\faq{\fa_\iq}
\def\faqp{\faq^+}
\def\faqd{\faq^*}
\def\faqcd{\fa_{\iq\siC}^*}

\def\Xp{X_+}
\def\oC{{}^\circ{\cC}}
\def\East{E^*}
\def\Ep{E_+}

\def\Parsigma{\Par_\sigma}

\def\ayF{{\ay_{F\iq}}}
\def\ayD{{\ay_{\Delta\iq}}}
\def\ayFd{{\ay_{F\iq}^*}}

\def\TFt{{\rmT}_F^t}

\def\TD{{\rmT}_\Delta}
\def\TDt{{\rmT}_\Delta^t}

\def\KDt{{\rmK}_\Delta^t}
\def\hinp#1#2{\langle#1|#2\rangle}
\def\Eastp{E^*_+}

\def\CuspD{\Cusp(X\col\tau)}
\def\CuspDt{\Cusp^t(X\col\tau)}
\def\CuspFv{\Cusp(\XFv\col\tau_F)}

\def\XFv{X_{F,v}}
\def\XFu{X_{F,u}}
\def\XFvp{X_{F,v,+}}

\def\Eps{E_{+,s}}

\def\oCFv{\oC_{F,v}}

\def\cWF{{}^F\cW}
\def\faqF{\ayF}
\def\faqFp{\fa^+_{F\iq}}

\def\baraqd{\bar{\fa}^*_\iq}

\def\stKF{\rmK_F^{\stt}}

\def\siC{{\scriptscriptstyle\complex}}

\def\siC{{\scriptstyle\complex}}
\def\faqdc{\fa_{\iq\siC}^*}
\def\faqFdperpc{\fa_{F\iq \siC}^{*\perp}}
\def\faqFdperp{\fa_{F\iq}^{*\perp}}

\def\faqFdc{\fa_{F\iq \siC}^*}
\def\faqFd{\fa_{F\iq}^*}

\def\PXt{\cP(X\col\tau)}
\def\PMXt{\cP_M(X\col\tau)}
\def\CMXt{C^\infty_M(X\col\tau)}
\def\Laufu{\Mer(\faqcd,\gS)^*_\laur}
\def\LauoC{\Laufu\otimes\oC^*}
\def\lsp{\bfL}
\def\ACX{\AC(X\col\tau)}
\def\ACR{\AC_\real(X\col\tau)}
\def\ACXFv{\AC(\XFv\col\tau_F)}
\def\PACX{\cP_{\AC}(X\col\tau)}
\def\PACC{\cP_{\AC}(X\col\tau)^\sim}
\def\cOS{\cO_S}
\def\cOSa{\cO_1}
\def\cOSo{\cO_0}
\def\cOSp{\cO'}
\def\CF{\Cusp_F}
\def\CFv{\CuspFv}
\def\nEFv{\nE_{F,v}}
\def\cTF{\cT^t_F}
\def\TstF{\rmT_F}

\def\PWXt{\PW(X\col\tau)}
\def\PWC{\PW(X\col\tau)^\sim}
\def\PWMX{\PW_M(X\col\tau)}
\def\Hypt{\Hyp^\sim}
\def\dt{d^\sim}
\def\cLa{\cL'_v}
\def\cLb{\cL}
\def\cLFv{\cL_v}
\def\reg{\mathop{\rm reg}}
\def\Hreg{H_{\reg}}
\def\ga{\alpha}

\def\geps{\varepsilon}  

\def\gl{\lambda}
\def\gL{\Lambda}
\def\gf{\varphi}

\def\gs{\sigma}
\def\gS{\Sigma}
\def\fa{\frak a}

\def\fg{\frak g}
\def\fh{\frak h}

\def\fk{\frak k}

\def\fn{\frak n}
\def\fp{\frak p}
\def\fq{\frak q}

\begin{document}
\currannalsline{164}{2006} 

\title{A Paley-Wiener theorem\\ for reductive symmetric spaces}

 \acknowledgements{}
\twoauthors{E. P. van den Ban}{H. Schlichtkrull}

 \institution{Mathematisch Instituut, Universiteit Utrecht, Utrecht, The
Netherlands\\
\email{ban@math.uu.nl}\\
\vglue-9pt
Matematisk Institut, K\o benhavns Universitet, K\o benhavn \O, Denmark
\\
\email{schlichtkrull@math.ku.dk}}

 \shorttitle{A Paley-Wiener theorem for reductive symmetric spaces} 
 
\centerline{\bf Abstract}
\vglue12pt

Let $X=G/H$ be a reductive symmetric space and $K$ a maximal
compact subgroup of $G$.
The image under the Fourier transform of the space of
$K$-finite compactly supported smooth functions on $X$
is characterized. 
\def\sni#1{\vskip1pt\noindent{#1}.\hskip5pt}
\vglue12pt
\centerline{\bf Contents} 
\vglue8pt
\sni{1} {Introduction}
\sni{2} {Notation}
\sni{3} {The Paley-Wiener space. Main theorem}
\sni{4} {Pseudo wave packets}
\sni{5} {Generalized Eisenstein integrals}
\sni{6} {Induction of Arthur-Campoli relations}
\sni{7} {A property of the Arthur-Campoli relations}
\sni{8} {Proof of Theorem 4.4}
\sni{9} {A comparison of two estimates}
\sni{10} {A different characterization of the Paley-Wiener space}

\section{Introduction}

One of the central theorems of harmonic analysis on $\real$ is the
Paley-Wiener theorem which characterizes the class of 
functions on $\complex$ which are Fourier transforms of
$C^\infty$-functions on $\real$ with compact support
(also called the Paley-Wiener-Schwartz theorem; see [18, p.\ 249]).
We consider the analogous question for
the Fourier transform of a reductive symmetric space $X=G/H$,
that is, $G$ is a real reductive Lie group of Harish-Chandra's class
and  $H$ is an open subgroup of the group $G^\gs$ of fixed points for 
an involution $\gs$ of $G.$ 

The paper is a continuation of [4] and [6],
in which we have shown that the Fourier transform is injective
on $\Cci(X)$, and established an inversion formula for the $K$-finite
functions in this space, with $K$ a $\sigma$-stable
maximal compact subgroup of $G$. 
A conjectural image of the space of $K$-finite functions in
$\Cci(X)$ was described in [4, Rem.~21.8], 
and will be confirmed 
in the present paper (the conjecture was already confirmed for symmetric 
spaces of split rank one in [4]).

If $G/H$ is a Riemannian symmetric space (equivalently, if $H$ is compact),
there is a well established theory of harmonic analysis (see
[17]), and the Paley-Wiener theorem that we obtain 
generalizes a well known theorem of Helgason and Gangolli
([15]; see also [17, Thm.\ IV,7.1]). 
Furthermore, the reductive group $G$ is a 
symmetric space in its own right, for the left times right action of 
$G \times G.$ Also in this `case of the group' there is an established 
theory of harmonic analysis, and our theorem generalizes the theorem 
of Arthur [1] (and Campoli [11] for groups of split rank one). 

The Fourier transform $\Fou$ that we are dealing with is defined for
functions in the space $\CciXt$ of $\tau$-spherical $\Cci$-functions 
on $X$. Here $\tau$ is a finite dimensional representation of
$K$, and a $\tau$-spherical function on $X$ is a function that
has values in the representation space $V_\tau$ and satisfies
$f(kx)=\tau(k)f(x)$ for all $x\in X$, $k\in K$. This space is
a convenient tool for the study of $K$-finite (scalar) functions
on $X$. Related to $\tau$ and the (minimal) principal series
for $X$, there is a family $\nE(\psi\col\gl)$
of normalized Eisenstein integrals on $X$ (cf.\ [2], [3]).
These are (normalized) generalizations of
the elementary spherical functions for Riemannian symmetric spaces,
as well as of Harish-Chandra's Eisenstein integrals
associated with a minimal parabolic subgroup of a semisimple Lie group.
The Eisenstein integral is a $\tau$-spherical smooth function on $X$.
It is linear in the parameter $\psi$, which
belongs to a finite dimensional Hilbert 
space~$\oC$, and meromorphic in 
$\gl$, which belongs to the complex linear dual $\faqcd$ of a maximal
abelian subspace $\faq$ of $\fp\cap\fq$. 
Here $\fp$ is the orthocomplement of $\fk$ in~$\fg$, and $\fq$ is the 
orthocomplement of $\fh$ in $\fg$,
where $\fg$, $\fk$ and $\fh$ are the Lie algebras of $G$, $K$ and $H$.
The Fourier transform $\Fou f$ of a function $f\in\CciXt$ is essentially
defined by integration of $f$ against $\nE$
(see (2.1)), and is a $\oC$-valued
meromorphic function of $\gl\in\faqcd$. The fact that $\Fou f(\lambda)$
is meromorphic in $\gl$, rather than holomorphic, represents a major 
complication not present in the mentioned special cases.

The Paley-Wiener theorem (Thm.\ 3.6) asserts that $\Fou$ maps $\CciXt$ 
onto the Paley-Wiener space $\PWXt$ (Def.\ 3.4), which is a space of 
meromorphic functions $\faqcd\to\oC$ 
characterized by an exponential growth condition and so-called
{\it Arthur-Campoli relations}, which are conditions coming from 
relations of a particular type among the Eisenstein integrals. These
relations generalize the relations used in [11] and [1].
Among the relations are conditions for
transformation under the Weyl group (Lemma 3.10). In the Riemannian case,
no other relations are needed, but this is not so in general.

The proof is based on the inversion formula $f=\cT\Fou f$ of [6], 
through which a function $f\in\CciXt$ is determined from its Fourier
transform by an operator $\cT$. The same operator can be applied to
an arbitrary function $\gf$ in the Paley-Wiener space $\PWXt$. The
resulting function $\cT\gf$ on $X$, called a pseudo wave packet, is
then shown to have $\varphi$ as its Fourier transform.
{\it A~priori\/}, $\cT\gf$ is defined and smooth on a certain dense open subset $X_+$ 
of $X$, and the main difficulty in the proof is to show that it
admits a smooth extension to $X$ (Thm.\ 4.4). 
In fact, as was shown already in
[6],   if a smooth extension of $\cT\gf$ exists, then
this extension has compact support and is mapped onto $\gf$ by $\Fou$.

The proof that $\cT\gf$ extends smoothly relies on the residue
calculus of [5] and on results of [7]. By means of the
residue calculus we write the pseudo wave packet
$\cT\gf$ in the form
$$\cT\gf=\sum_{F\subset\Delta} \cT_F\gf$$(see eq.\ (8.3))
in which $\Delta$ is a set of simple roots for the root system of
$\faq$, and in which
the individual terms  for $F\neq\emptyset$
are defined by means of residue operators. The term
$\cT_\emptyset\gf$ is the wave packet given by integration over $\faqd$
of $\gf$ against the normalized Eisenstein integral. The smooth
extension of $\cT\gf$ is established by showing that each term
$\cT_F\gf$ extends smoothly. The latter fact is obtained by
identification of $\cT_F\gf$ with a wave packet formed by
generalized Eisenstein integrals. The generalized Eisenstein 
integrals we use were introduced in [6]; they are 
smooth functions on $X$. It is shown in [9]
that they are matrix coefficients of nonminimal principal
series representations and that they agree with the generalized
Eisenstein integrals of [12].
However, these facts play no role here.
It is for the identification of $\cT_F\gf$ 
as a wave packet that the Arthur-Campoli
relations are needed when $F\neq\emptyset$. 
An important step is to show that
Arthur-Campoli relations for lower dimensional symmetric spaces,
related to certain parabolic subgroups in $G$, can be induced up
to Arthur-Campoli relations for $X$ (Thm.\ 6.2). For this
step we use a result from [7].

As mentioned, our Paley-Wiener theorem generalizes 
that of Arthur [1] for the group case. 
Arthur also uses residue calculus in the spirit of [19],
but apart from that our approach differs
in a number of ways, the following two being the most
significant.  
Firstly, Arthur relies on Harish-Chandra's Plancherel theorem for
the group, whereas we do not need the analogous
theorem for $X$, which has been established by Delorme [14] and
the authors [8], [9]. 
Secondly, Arthur's result involves unnormalized 
Eisenstein integrals, whereas our involves normalized ones.
This facilitates comparison between the Eisenstein integrals 
related to $X$ and those related to lower rank symmetric spaces coming
from parabolic subgroups. 
For similar comparison of the unnormalized Eisenstein integrals,
Arthur relies on a lifting principle of Casselman, the proof of which
has not been published. In [7] we have established a normalized version
of Casselman's principle which plays a crucial role in the present work.
One can show, using [16, Lemma 2, p.\ 156], [1, Lemma I.5.1]
and [13], that our Paley-Wiener
theorem, specialized to the group case, implies Arthur's. In
fact, it implies a slightly stronger result, since here only
Arthur-Campoli relations for real-valued parameters $\gl$ are needed,
whereas the Paley-Wiener theorem of [1] requires also
the relations at the complex-valued $\gl$.

The Paley-Wiener space $\PWXt$ is defined in Section 3
(Definition 3.4), and the
proof outlined above that it equals the 
Fourier image of $\CciXt$ takes up the following Sections 4--8.
{\it A priori\/} the given definition of $\PWXt$ does not match that of
[4], but it is shown in the final Sections 9, 10  that
the two spaces are equal.

The main result of this paper was found and announced in the fall of 1995
when both authors were visitors of the Mittag-Leffler Institute in Djursholm, 
Sweden. We are grateful to the organizers of the program and the staff of 
the institute for providing us with this opportunity, and to Mogens 
Flensted-Jensen for helpful discussions during that period.

\vglue-18pt
\phantom{up}
\section{Notation}
\vglue-4pt

We use the same notation and basic assumptions as
in [4, \S\S 2, 3, 5, 6], and [6, \S 2]. 
Only the most essential notions will be recalled, and we refer to
the mentioned locations for unexplained notation.

We denote by $\gS$ the root system of $\faq$ in $\fg$, where $\faq$
is a maximal abelian subspace of $\fp\cap\fq$,
as mentioned in the introduction. 
Each positive system $\gS^+$ for $\gS$ determines a 
parabolic subgroup $P=M_1N$, where $M_1$ is the centralizer 
of $\faq$ in $G$ and $N$ is the exponential of $\fn$, the sum
of the positive root spaces. In what follows we
assume that such a positive system $\gS^+$ has been fixed.
Moreover, notation  with reference to $\gS^+$ or $P$, as given
in [4] and [6], is 
supposed to refer to this fixed choice, if nothing else is mentioned.
For example, we write $\faqp$  for the corresponding
positive open Weyl chamber in $\faq$, denoted $\faqp(P)$  in
[4], and $\Aqp$ for its exponential $\Aqp(P)$ in $G$. 
We write $P=MAN$ for the Langlands decomposition of~$P$.

Throughout the paper we fix a finite dimensional unitary 
representation $(\tau,\Vtau)$ of $K$, and we denote by $\oC=\oC(\tau)$ 
the finite dimensional space defined by [4, eq.\ (5.1)].
The Eisenstein integral $E(\psi\col\gl)=E(P\col\psi\col\gl)\colon X\to\Vtau$ 
is defined as in [4, eq.\ (5.4)], and the normalized Eisenstein 
integral $\nE(\psi\col\gl)=\nE(P\col\psi\col\gl)$ is defined as in 
[4, p.\ 283]. 
Both Eisenstein integrals belong to $\Cinf(X\col\tau)$ and depend
linearly on $\psi\in\oC$ and meromorphically on $\gl\in\faqdc$. 
For $x\in X$ we denote
the linear map $\oC\ni\psi\mapsto\nE(\psi\col\gl\col x)\in\Vtau$
by $\nE(\gl\col x)$, and we define $\East(\gl\col x)\in\Hom(\Vtau,\oC)$
to be the adjoint of $\nE(-\bar\gl\col x)$ (see [6, eq.\ (2.3)]).
The Fourier transform that we investigate  maps
$f\in\CciXt$ to the meromorphic function $\Fou f$ on $\faqdc$ given by
\begin{equation}
\Fou f(\gl)= 
\int_X \East(\gl\col x) f(x)\,dx\in \oC.
\end{equation}

The open dense set $\Xp\subset X$ 
is given by
$$\Xp=\cup_{w\in\cW}\, K\Aqp wH;$$
see [6, eq.\ (2.1)]. It naturally arises in connection with
the study of asymptotic expansions of the Eisenstein integrals;
see [6, p.\ 32, 33]. As a result of this theory, the normalized
Eisenstein integral is decomposed as a finite sum
\begin{equation}
\nE(\gl\col x)=\sum_{s\in W}\Eps(\gl\col x), \qquad
\Eps(\gl\col x)=\Ep(s\gl\col x)\after\nC(s\col\gl)
\end{equation}
for $x\in\Xp$, all ingredients being
meromorphic in $\gl\in\faqdc$. The partial Eisenstein integral 
$\Ep(\gl\col x)$ is a 
$\Hom(\oC,V_\tau)$-valued function in $x\in\Xp$, given 
by a converging series expansion, 
and $\nC(s\col\gl)\in\End(\oC)$ is the (normalized)\break
$c$-function associated with $\tau$. In general,
$x\mapsto\Ep(\gl\col x)$ is singular along $X\setminus \Xp$.
The $c$-function also appears in the following transformation law
for the action of the Weyl group
\begin{equation}\East(s\gl\col x)=\nC(s\col\gl)\after\East(\gl\col x)
\end{equation}
 for all
$s\in W$ and $x\in X$ (see [6, eq.\ (2.11)]), from which it follows
that 
\begin{equation}\Fou f(s\gl)=\nC(s\col\gl)\after\Fou f(\gl).
\end{equation}

The structure of the singular set for the meromorphic functions
$\nE(\,\cdot\,\col x)$ and $\Ep(\,\cdot\,\col x)$ on $\faqcd$
plays a crucial role. To describe it, we
recall from [7, \S 10], that a $\gS$-configuration 
in $\faqcd$ is a locally finite collection of affine hyperplanes
$H$ of the form 
\begin{equation}
H=\{\gl\mid
\langle\gl,\ga_H\rangle=s_H\}
\end{equation}
 where
$\ga_H\in\gS$
and $s_H\in\complex$. 
Furthermore, we recall from [7, \S 11], that
if $\Hyp$ is a $\gS$-configuration 
in $\faqcd$ and $d$ a map $\Hyp\to\Natural$, we define
for each bounded set $\omega\subset\faqcd$ a polynomial function
$\pi_{\omega,d}$ on $\faqcd$ by 
\begin{equation}
\pi_{\omega,d}(\gl)=\prod_{H\in\Hyp, H\cap\omega\neq\emptyset}
(\langle\gl,\ga_H\rangle-s_H)^{d(H)},
\end{equation}
where $\ga_H, s_H$ are as above.
The linear space $\Mer(\faqcd,\Hyp,d)$ 
is defined to be the space of meromorphic functions
$\gf\colon\faqcd\to\complex$, for which $\pi_{\omega,d}\gf$
is holomorphic on $\omega$ for all
bounded open sets $\omega\subset\faqcd$, and
the linear space $\Mer(\faqcd,\Hyp)$ is defined
by taking the union of $\Mer(\faqcd,\Hyp,d)$ over $d\in\Natural^\Hyp$.
If $\Hyp$ is real, that is, 
$s_H\in\real$ for all $H$, we write
$\Mer(\faqd,\Hyp,d)$ and $\Mer(\faqd,\Hyp)$
in place of $\Mer(\faqcd,\Hyp,d)$ and $\Mer(\faqcd,\Hyp)$.

\begin{lem}
There exists a real $\gS$-configuration $\Hyp$ such that
the meromorphic functions
$\nE(\,\cdot\,\col x)$
and
$\Eps(\,\cdot\,\col x')$ belong to $\Mer(\faqd,\Hyp)\otimes\Hom(\oC,\Vtau)$
for all $x\in X${\rm ,} $x'\in X_+${\rm ,} $s\in W${\rm ,} and such that
$\nC(s\col\,\cdot\,)\in\Mer(\faqd,\Hyp)\otimes\End(\oC)$ 
for all $s\in W$.
\end{lem}  

\Proof The statement for $\nE(\,\cdot\,\col x)$ is proved in
[6, Prop.\ 3.1], and the statement for
$E_{+,1}(\,\cdot\,\col x)=
\Ep(\,\cdot\,\col x)$ is proved in
[6, Lemma 3.3]. The statement about  $\nC(s\col\,\cdot\,)$
follows from [3, eqs.\ (68), (57)],
by the argument given below the proof of Lemma 3.2 in
[6]. The statement for
$\Eps(\,\cdot\,\col x)$ in general then follows from its definition
in (2.2). 
\Endproof 

Let $\Hyp=\Hyp(X,\tau)$ denote the collection of the singular
hyperplanes for all $\gl\mapsto\East(\gl\col x)$, $x\in X$ 
(this is a real $\gS$-configuration, by the preceding lemma). 
Moreover, for $H\in\Hyp$ let $d(H)=d_{X,\tau}(H)$ be the least 
integer $l\ge0$ for which
$\gl\mapsto(\langle\gl,\ga_H\rangle-s_H)^{l}\East(\gl\col x)$ is regular
along $H\setminus\cup\{H'\in\Hyp\mid H'\neq H\}$, for all $x\in X$.
Then
$\East(\,\cdot\,\col x)\in\Mer(\faqd,\Hyp,d)\otimes\Hom(\Vtau,\oC)$
and $d$ is minimal with this property. It follows that
$\Fou f\in\Mer(\faqd,\Hyp,d)\otimes\oC$ for all $f\in\CciXt$.

There is more to say about these singular sets.
For $R\in\real$ we define
\begin{equation}
\faqd(P,R)=\{\gl\in\faqcd\mid\forall\alpha\in\gS^+:\,
\Re\langle\gl,\alpha\rangle<R
\}
\end{equation}
and denote by $\baraqd(P,R)$ the closure of this set.
Then it also follows from\break [6, Prop.\ 3.1 and  Lemma 3.3],
that $\East(\,\cdot\,\col x)$ and $\Ep(\,\cdot\,\col x)$
both have the property that for each $R$
only finitely many singular hyperplanes 
meet $\faqd(P,R)$.

In particular, the set
of affine hyperplanes
\begin{equation}
\Hyp_0=\{H\in\Hyp(X,\tau)\mid H\cap\baraqd(P,0)\neq\emptyset\},
\end{equation}
is finite.
Let $\pi$ be the real polynomial function on $\faqcd$ given by
\begin{equation}
\pi(\gl)=\prod_{H\in\Hyp_0}
(\langle\gl,\ga_H\rangle-s_H)^{d_{X,\tau}(H)}
\end{equation}
where $\ga_H$ and $s_H$ are chosen as in (2.5).
The polynomial $\pi$
coincides, up to a constant nonzero factor, with the polynomial 
denoted by the same symbol in [4, eq.\ (8.1)], and in [6,
p.\ 34]. 
It has the
property that there exists $\geps>0$ such that $\gl\mapsto\pi(\gl)
\East(\gl\col x)$ is holomorphic on 
$\faqd(P,\geps)$ 
for all $x\in X$.

\vglue-20pt
\phantom{up}
\section{The Paley-Wiener space. Main theorem}
\vskip-6pt

We define the Paley-Wiener space
$\PWXt$ for the pair $(X,\tau)$ and state the main theorem, 
that the Fourier transform maps $\CciXt$ onto
this space.

First we set up the condition that reflects relations among Eisenstein
integrals.
In [11] and [1] similar relations are used in the definition
of the Paley-Wiener space.
However, as we are dealing with functions that are in general
meromorphic rather than holomorphic, 
our relations have to be specified somewhat differently.
This is done by means of Laurent functionals,
a concept introduced in [7, Def.\ 10.8], to which we refer
(see also the review in [8, \S 4]).
In [4, Def.\ 21.6], the required relations are
formulated differently; we compare the definitions in Lemma
10.4 below.

\vglue-19pt
\phantom{up}

\begin{defi} We call a $\gS$-Laurent functional 
$\cL\in\LauoC$ an {\it Arthur-Campoli
functional} if it annihilates $\East(\,\cdot\,\col x)v$ for all $x\in X$
and $v\in V_\tau$. The set of all Arthur-Campoli
functionals is denoted
$\ACX$, and the subset of the Arthur-Campoli functionals with support
in $\faqd$ is denoted
$\ACR$.\end{defi}

It will be shown below in Lemma 3.8 that the elements
of $\ACX$ are natural objects, 
from the point of view of characterizing $\Fou(\CciXt)$.

Let $\Hyp$ be a real $\gS$-configuration in $\faqdc$, and
let $d\in\Natural^\Hyp$.
By $\Pal(\faqd,\Hyp,d)$ we denote the linear space of functions
$\gf\in\Mer(\faqd,\Hyp,d)$ with polynomial decay in the
imaginary directions, that is
\begin{equation}\sup_{\gl\in\omega+i\faqd} (1+|\gl|)^n|\pi_{\omega,d}(\gl)\gf(\gl)|
<\infty
\end{equation}
for all compact $\omega\subset\faqd$ and all $n\in\Natural$. 
The space $\Pal(\faqd,\Hyp,d)$ is given
a Fr\'echet topology by means of the seminorms in 
(3.1).
The union of these spaces over all $d\colon\Hyp\to\Natural$,
equipped with the limit topology,
is denoted $\Pal(\faqd,\Hyp)$.

\begin{defi} Let $\Hyp=\Hyp(X,\tau)$ and $d=d_{X,\tau}$.
We define 
$$\PACX=\{\gf\in\Pal(\faqd,\Hyp,d)\otimes\oC\mid\cL\gf=0, 
\forall \cL\in\ACR\},$$
and equip this subspace of $\Pal(\faqd,\Hyp,d)\otimes\oC$ with the 
inherited topology.
\end{defi}

\begin{lem}
The space $\PACX$ is a Fr\ee chet space.\end{lem}

\Proof  Indeed, $\PACX$ is a closed subspace of 
$\Pal(\faqd,\Hyp,d)\otimes\oC$, since Laurent functionals are continuous
on $\Pal(\faqd,\Hyp,d)$ (cf.\ [5, Lemma 1.11]).\Endproof\vskip4pt 

In Definition 3.2 it is required that the elements of $\PACX$ belong to
$\Pal(\faqd,\Hyp,d)\otimes\oC$ where $\Hyp=\Hyp(X,\tau)$ and 
$d=d_{X,\tau}$ are specifically given in terms of the singularities
of the Eisenstein integrals.
It will be shown in Lemma 3.11 below
that this requirement 
is unnecessarily strong (however, it is convenient for the definition
of the topology).

\begin{defi}
The {\it Paley-Wiener space\/} $\PWXt$ is defined
as the space of functions $\gf\in\PACX$ for which there
exists a constant $M>0$ such that
\begin{equation}
\sup_{\gl\in\baraqd(P,0)} (1+|\gl|)^n e^{-M\,|\!\Re\gl|}
\|\pi(\gl)\gf(\gl)\|<\infty
\end{equation}
for all $n\in\Natural$. The subspace of functions
that satisfy (3.2) for all $n$ and
a fixed $M>0$ is denoted $\PWMX$. The space $\PWMX$ is given
the relative topology of $\PACX$, or equivalently,
of $\Pal(\faqd,\Hyp,d)\otimes\oC$ 
where $\Hyp=\Hyp(X,\tau)$ and $d=d_{X,\tau}$. 
Finally, the Paley-Wiener space $\PWXt$ is given the
limit topology of the union 
\begin{equation}
\PWXt=\cup_{M>0}\PWMX.
\end{equation}
\end{defi}

The functions in $\PWXt$ are called {\it Paley-Wiener functions}.
By the definition just given they are the functions in
$\Mer(\faqd,\Hyp,d)\otimes\oC$ for which the estimates (3.1) and
(3.2) hold, and which are annihilated by all Arthur-Campoli functionals
with real support.

\begin{rem} It will be verified later that $\PWMX$ is a closed
subspace of $\PACX$ (see Remark 4.2).
Hence $\PWMX$ is a Fr\'echet space, and $\PWXt$ a strict LF-space
(see [20, p.~291]). Notice that the Paley-Wiener space $\PWXt$
is not given the relative topology of $\PACX$. However, the
inclusion map $\PWXt\to\PACX$ is continuous.
\end{rem}

We are now able to state the 
Paley-Wiener theorem for the pair $(X,\tau)$. 

\begin{theorem} The Fourier transform $\Fou$ is
a topological linear isomorphism of $\CMXt$ onto 
$\PWMX${\rm ,} for each $M>0${\rm ,}
and it is a topological linear isomorphism of $\CciXt$ onto 
the Paley-Wiener space $\PWXt$.\end{theorem}

Here we recall from [6, p.\ 36], that $\CMXt$ is the subspace of
$\Cinf(X\col\tau)$ consisting of those functions that are supported
on the compact set $K\exp B_MH$, where $B_M\subset\faq$ is the closed
ball of radius $M$, centered at 0. The space $\CMXt$
is equipped with its standard Fr\'echet topology, which
is the relative topology of
$\Cinf(X\col\tau)$. Then
\begin{equation}
\Cci(X\col\tau)=\cup_{M>0}\CMXt 
\end{equation}
and $\Cci(X\col\tau)$ carries the limit topology of this union.

The final statement in the theorem is an obvious consequence of the first,
in view of (3.3) and (3.4). The proof of the first statement 
will be given in the course of the next 5 sections 
(Theorems 4.4, 4.5, proof in Section~8).
It relies on several results from [6], which are elaborated
in the following two sections. At present, we note the following:

\begin{lem} The Fourier transform $\Fou$ maps $\CMXt$
continuously and injectively into $\PWMX$ for each $M>0$.\end{lem}

\Proof The injectivity of $\Fou$ is one of the main results in
[4,   Thm.\ 15.1].
It follows from [6, Lemma 4.4], that $\Fou$ maps
$\CMXt$ continuously into the space $\Pal(\faqd,\Hyp,d)\otimes\oC$,
where $\Hyp=\Hyp(X,\tau)$ and $d=d_{X,\tau}$,
and that (3.2) holds for $\gf=\Fou f\in\Fou (\CMXt)$.
Finally, it follows from Lemma 3.8 below that $\Fou$ maps into $\PACX$.
\hfill\qed

\begin{lem} Let $\cL\in\LauoC$. Then $\cL\in\ACX$ if and only if
$\cL\Fou f=0$ for all $f\in\CciXt$.
\end{lem}

\Proof Recall that $\Fou f$ is defined by (2.1)
for $f\in\CciXt$. We claim that 
\begin{equation}\cL\Fou f= 
\int_X \cL\East(\,\cdot\,\col x) f(x)\,dx,
\end{equation}
that is, the application of $\cL$ can be taken inside the integral.

The function $\lambda\mapsto \East(\lambda\col x)$
on $\faqcd$ belongs to $\Mer(\faqd,\cH,d)\otimes\oC$ for each
$x\in X$, where $\Hyp=\Hyp(X,\tau)$ and $d=d_{X,\tau}$.
The space $\Mer(\faqd,\cH,d)\otimes\oC$ is a complete locally
convex space, when equipped with the initial topology with respect to
the family of maps $\gf\mapsto\pi_{\omega,d}\gf$ into $\cO(\omega)$,
and $x\mapsto\East(\,\cdot\,\col x)$
is continuous (see [3, Lemma 14]). The integrals in
(2.1) and (3.5) may be seen as integrals with values in this
space. Since Laurent functionals
are continuous, (3.5) is justified.

Assume now that $\cL\in\ACX$ and let $f\in\CciXt$. Then\break
$\cL\East(\,\cdot\,\col x)f(x)=0$ for each $x\in X$, and the 
vanishing of $\cL\Fou f$ follows immediately
from (3.5).

Conversely, assume that $\cL$ annihilates $\Fou f$ for all $f\in\CciXt$.
From (3.5) and
[4, Lemma 7.1], it follows easily that $\cL$ annihilates 
$\East(\,\cdot\,\col a)v$ for\break\vskip-12pt\noindent $v\in V_\tau^{K\cap H\cap M}$ and
$a\in\Aqp(Q)$, with $Q\in\Parmin$ arbitrary. Let $v\in V_\tau$. Since
$\East(\gl\col kah)= \East(\gl\col a)\after\tau(k)^{-1}$
for $k\in K$, $a\in\Aq$ and $h\in H$, it is seen that
$\East(\gl\col kah)v=\East(\gl\col a)P(\tau(k)^{-1}v)$
where $P$ denotes the orthogonal projection $V_\tau\to
V_\tau^{K\cap H\cap M}.$ Hence $\cL$ 
annihilates $\East(\,\cdot\,\col x)v$ for all $x\in\Xp$, $v\in V$.
By continuity and density 
the same conclusion holds for all $x\in X$.
\hfq

\begin{rem} In Definition 3.2 we used only Arthur-Campoli
functionals with real support. Let $\PACC$ denote the space obtained
in that definition with $\ACR$ replaced by $\ACX$, and let $\PWC$
denote the space obtained in Definition 3.4{} with $\PACX$
replaced by $\PACC$. Then clearly $\PACC\subset\PACX$ and
$\PWC\subset\PWXt$. However, it follows from Lemma 3.8
that $\Fou(\CciXt)\subset\PWC$, and hence as a consequence of
Theorem 3.6{} we will have
$$\PWC=\PWXt.$$
\end{rem}

In general, the Arthur-Campoli functionals are not explicitly 
described. Some relations of a more explicit nature can be pointed 
out: these are the relations (2.4) that express 
transformations
under the Weyl group. In the following lemma it is shown
that these relations are of Arthur-Campoli type, which explains
why they are not mentioned separately in the definition of
the Paley-Wiener space. 

\begin{lem} Let $\gf\in\PACX$. Then
$\gf(s\gl)=\nC(s\col\gl)\gf(\gl)$
for all $s\in W$ and $\gl\in\faqcd$ generic.\end{lem}

\Proof The relation 
$\gf(s\gl)=\nC(s\col\gl)\gf(\gl)$
is meromorphic in $\gl$, so it suffices to verify it for
$\gl\in\faqd$.
Let $\Hyp=\Hyp(X,\tau)$.
Fix $s\in W$ and $\gl\in\faqd$ such that $\nC(s\col\gl)$ is
nonsingular at $\gl$, and such that $\gl$ and $s\gl$ 
do not belong to any of the hyperplanes from $\Hyp$. 
Let $\psi\in\oC$ and consider the linear
form $\cL_\psi\maps\gf\mapsto\hinp{\gf(s\gl)-\nC(s\col\gl)\gf(\gl)}{\psi}$ on
$\Mer(\faqd,\Hyp)\otimes\oC$. It follows from [7, Remark 10.6],
that for each $\nu\in\faqdc$ there exists
a $\gS$-Laurent functional which, when applied to the functions that
are regular at $\nu$, yields the evaluation in $\nu$. Obviously,
the support of such a functional is $\{\nu\}$.
Hence there exists $\cL\in\LauoC$ with support
$\{\gl,s\gl\}$ such that 
$\cL\gf=\cL_\psi\gf$ for all 
$\gf\in\Mer(\faqd,\Hyp)\otimes\oC$.
It follows from (2.3) and
Definition~3.1  that $\cL\in\ACR$. The lemma follows immediately.
\hfq

\begin{lem} Let $\Hyp$ be a real $\gS$-configuration
in $\faqcd$ and let $\gf\in\Pal(\faqd,\Hyp)\break\otimes\oC$.
Assume $\cL\gf=0$ for all $\cL\in\ACR$. Then $\gf\in\PACX$.\end{lem}

\Proof Let $d\in\Natural^\Hyp$ be such that 
$\gf\in\Pal(\faqd,\Hyp,d)\otimes\oC$.
We may assume that $\Hyp\supset\Hyp(X,\tau)$ and that
$d\succeq d_{X,\tau}$ (that is,
$d(H)\geq d_{X,\tau}(H)$ for all $H\in\Hyp$), 
where $d_{X,\tau}$ is trivially
extended to $\Hyp$. Let $H\in\Hyp$ be arbitrary and
let $l$ be the least nonnegative integer for which 
$\gl\mapsto(\langle\gl,\ga_H\rangle-s_H)^{l}\gf(\gl)$
is regular along  $\Hreg\coleq H\setminus\cup\{H'\in\Hyp\mid H'\neq H\}$. 
Then $l\le d(H)$, and 
the statement of the lemma amounts to
$l\le d_{X,\tau}(H)$. 

Assume that $l>d_{X,\tau}(H)$; we will show that this leads
to a contradiction.
Let $d'\in\Natural^\Hyp$ be the
element such that $d'(H)=l$ and which equals $d$ on all other hyperplanes
in $\Hyp$. Then $\gf\in\Pal(\faqd,\Hyp,d')\otimes\oC$ and
$d'\succ d_{X,\tau}$. 
Let $\gl_0\in\Hreg\cap\faqd$.
It follows from [7, Lemmas 10.4, 10.5],
that there exists 
$\cL\in\Laufu$ such that $\cL\phi$ is the evaluation in $\gl_0$ of
$(\langle\gl,\ga_H\rangle-s_H)^{l}\phi(\gl)$ 
for all $\phi\in\Mer(\faqd,\Hyp, d')$.
Obviously, $\supp\cL=\{\gl_0\}\subset\faqd$. 
Since $l>d_{X,\tau}(H)$,
the functional $\cL\otimes\eta$ annihilates
$\Mer(\faqd,\Hyp, d_{X,\tau})\otimes\oC$ for all $\eta\in\oC^*$ and hence
belongs to $\ACR$. Then it also annihilates $\gf$, that is, 
the function $(\langle\gl,\ga_H\rangle-s_H)^{l}\gf(\gl)$ vanishes
at 
$\gl_0$, which was arbitrary in $\Hreg\cap\faqd$. 
By meromorphic continuation this function 
vanishes everywhere. 
This contradicts the definition of~$l$.
\hfq

\section{Pseudo wave packets}

In the Fourier inversion formula $\psW\Fou f=f$
the 
{\it pseudo wave packet} $\psW\Fou f$ is defined by
\begin{equation}
\psW\Fou f(x)=|W|\int_{\eta+i\faqd} \Ep(\gl\col x)\Fou f(\gl)\,d\gl,\qquad
x \in\Xp,  
\end{equation}
for $f\in\CciXt$ and for $\eta\in\faqd$ sufficiently antidominant (the function is then 
independent of $\eta$). 
Here $d\gl$ is the translate of Lebesgue measure on $i\faqd$,
normalized as in [6, eq.\ (5.2)].
{\it A priori\/}, $\psW\Fou f$ belongs
to the space
$\Cinf(\Xp\col\tau)$ of smooth $\tau$-spherical functions on $\Xp$, but
the identity with $f$ shows that it extends
to a smooth function on $X$.

The pseudo wave
packets are also used for the proof of the Paley-Wiener theorem: 
Given a function in the Paley-Wiener space, the candidate for
its Fourier preimage is constructed as a pseudo wave packet on $\Xp$.
In this section we reduce the proof of the Paley-Wiener theorem 
to one property
of such pseudo wave packets.
This property, that they extend to global smooth 
functions on $X$, will be established in Section~8

We first recall some spaces defined in [6], and relate
them to the spaces given in Definitions~3.2 and 3.4.

\begin{defi}
Let $\PXt$ be the space of meromorphic functions $\gf\maps\faqcd\break\to\oC$
having the following properties (i)--(iii) (see (2.9) for the definition of 
$\pi$):
\begin{itemize}
\ritem{(i)} $\gf(s\gl)=\nC(s\col\gl)\gf(\gl)$ for all $s\in W$
and generic $\gl\in\faqdc$.

\ritem{(ii)} There exists $\geps>0$ such that $\pi\gf$ is
holomorphic on $\faqd(P,\geps)$.

\ritem{(iii)} For some $\epsilon>0$, for every compact
set $\omega\subset\faqd(P,\geps)\cap\faqd$ and for all $n\in\Natural$,
$$\sup_{\gl\in\omega+i\faqd} (1+|\gl|)^n \|\pi(\gl)\gf(\gl)\|<\infty.$$
\end{itemize}
\noindent Moreover, for each $M>0$ let
$\PMXt$ be the subspace of $\PXt$ consisting of the functions
$\gf\in\PXt$ with the following property (iv). 
\begin{itemize}
\ritem{(iv)} For every
strictly antidominant $\eta\in\faqd$ there exists a constant
$t_\eta\ge 0$ such that\end{itemize}
\begin{equation}
\sup_{t\ge t_\eta, \gl\in t\eta+i\faqd}
(1+|\gl|)^{\dim\faq+1}e^{-M\,|\!\Re\gl|}\|\gf(\gl)\|<\infty.
\end{equation}
\end{defi}

Notice that (ii) and (iii) are satisfied by any function $$\gf\in
\Pal(\faqd,\Hyp(X,\tau),d_{X,\tau})\otimes\oC,$$ by the definition of $\pi$. 
If $\gf$ belongs to the subspace
$\PACX$ it also satisfies (i), by Lemma 3.10, and hence
\begin{equation}
\PWXt\subset\PACX\subset\PXt.
\end{equation}
Moreover, the estimate in
(3.2) is stronger than (iv), and hence 
\begin{equation}
\PWMX\subset\PACX\cap\PMXt.
\end{equation}

\begin{rem}
 It will be shown later by Euclidean Fourier analysis,
see Lemma~9.3, that the stronger estimate
(3.2) holds for all $\gf\in\PMXt$. In
particular, it follows that in fact
\begin{equation}
\PWMX=\PACX\cap\PMXt.
\end{equation}
It will also follow from Lemma~9.3  that 
$\PWMX$ is a closed subspace of $\PACX$, hence
a Fr\'echet space. Alternatively, the latter property of $\PWMX$ follows 
directly from Theorem 3.6, in the proof of which it is never
used. In fact, (4.5) will be established
in the course of that proof. 
\end{rem}

\begin{rem} It will also be shown, see Lemma 10.2, that there exist 
a real $\gS$-configura\-tion $\Hypt$ and a map 
$\dt\colon\Hypt\to\Natural$ such that 
$\PXt\subset\Pal(\faqd,\Hypt,\dt)\otimes\oC$. In combination with
Lemma 3.11 this implies that
$$\PACX=\{\gf\in\PXt\mid \cL\gf=0,\forall\cL\in\ACR\}.$$
The 
present remark is not used in the proof of Theorem~3.6. 
\end{rem}

Recall from [6, \S 4], that the 
{\it pseudo wave packet} of (4.1) can be formed with
$\Fou f$ replaced by an arbitrary 
function $\gf\in\PXt$. The resulting function
$\psW\gf\in\Cinf(\Xp\col\tau)$ is given by 
\begin{equation}
\psW\gf(x)=|W|\int_{\eta+i\faqd} \Ep(\gl\col x)\gf(x)\,d\gl,\qquad
x \in\Xp,
\end{equation}
for $\eta\in\faqd$ sufficiently antidominant, so that
the function is independent of~$\eta$. The following theorem
represents the main step in the proof of the Paley-Wiener theorem.

\begin{theorem} Let $\gf\in\PACX$.
Then $\psW\gf$ 
extends to a smooth\break $\tau$-spherical function on $X$
\/{\rm (}\/also denoted by $\psW\gf${\rm ). }
The map $\psW$ is continuous from
$\PACX$ to $\Cinf(X\col\tau)$.\end{theorem}

We will prove this result in Section 8{} (see below Theorem 8.3).
However, we first use it to derive the following Theorem 4.5, from
which Theorem 3.6{} is an immediate consequence. 

\begin{theorem} \hskip-5pt  Let $M\!>\!0$. Then $\psW\gf\!\in\!\CMXt$ for all $\gf\!\in\!\PWMX${\rm ,}
and $\psW$ 
is a continuous inverse to the Fourier transform $\Fou\maps\CMXt\to\PWMX$.
\end{theorem} 

\Proof Let $\cP'_M(X\col\tau)$ denote the set of functions 
$\gf\in\PMXt$ for which $\psW\gf$ has a smooth extension to $X$.
We have seen in [6, Cor.\ 4.11], that $\Fou$ maps $\CMXt$ bijectively
onto $\cP'_M(X\col\tau)$ with $\psW$ as its inverse. It follows from
Theorem 4.4{} that $\PACX\cap\PMXt$ is contained in $\cP'_M(X\col\tau)$.
Combining this with Lemma 3.7{} and (4.4) we obtain the following
chain of inclusions
\begin{eqnarray*}
\Fou(\CMXt)&\subset&\PWMX\subset\PACX\cap\PMXt\\
&\subset&\cP'_M(X\col\tau)
=\Fou(\CMXt).
\end{eqnarray*}
It follows that these inclusions are equalities
(in particular, (4.5) is then established). Thus $\Fou$ is
bijective  $\CMXt\to\PWMX$, with inverse $\psW$.

Since $\psW\maps\PACX\to\Cinf(X\col\tau)$ is continuous by Theorem
4.4  and since $\PWMX$ and $\CMXt$ carry the restriction 
topologies of these spaces,
we conclude that the restriction map $\psW\maps\PWMX\to\CMXt$ is continuous.
\phantom{endofline}\hfq

\section{Generalized Eisenstein integrals}

In [6, \S 10], we defined generalized Eisenstein integrals
for $X$. These will be used extensively in the following. In this
section we recall their definition and derive some properties of them.
For further properties (not to be used here), we refer to
[8], [9].

Let $t\in\WT(\gS)$ be an even and $W$-invariant residue weight
(see [5, p.\ 60]) to be fixed throughout the paper.
Let
$f\mapsto \TDt f$, $\CciXt\to\Cinf(X\col\tau)$, be the operator
defined by [6, eq.\ (5.5)], with $F=\Delta$. The fact that it maps into 
$\Cinf(X\col\tau)$ is a consequence of [6, Cor.\ 10.11].
Moreover, if the vectorial part of $X$ vanishes, that is, if
$\ayD=\{0\}$, then
\begin{equation}
\TDt f(x)=|W|
\int_X \KDt(x \col y) f(y)\,dy
\end{equation}
for $x\in X$, cf.\ [6, eq.\ (5.10) and proof of Cor.\ 10.11],
where $\KDt(x \col y)$ is the residue kernel defined by
[6, eq.\ (5.7)], with $F=\Delta$.

If the vectorial part of $X$ vanishes, 
then we follow [6, Remark 10.5], and
define a finite dimensional space by 
\begin{equation}
\CuspDt=\Span\{\KDt(\,\cdot\,\col y)u\mid
y\in\Xp,u\in\Vtau\}\subset\Cinf(X\col\tau).
\end{equation}
The space is denoted $\cC_\Delta$ in [6], whereas the present
notation is in agreement with [8, \S 9].
By continuity of $\KDt$ and finite dimensionality of $\CuspDt$, 
$\KDt(\,\cdot\,\col y)u$ belongs to this space for 
$y\in X\setminus\Xp$  as well.

\begin{lem}
Assume $\ayD=\{0\}$. Then $\TDt f\in\CuspDt$ 
for all $f\in\CciXt${\rm ,} and 
the map $\TDt\colon\CciXt\to\CuspDt$ is surjective.
\end{lem}

\Proof The map $y\mapsto \KDt(\,\cdot\,\col y)f(y)$ belongs to
$\CciXt\otimes\CuspDt$. Hence its integral (5.1) over $X$ belongs 
to $\CuspDt$. The surjectivity follows from (5.1); see 
[8, Lemma 9.1].
\hfq

\begin{rem} It is seen in [8, Thm.\ 21.2,
Def.\ 12.1 and Lemma 12.6],
that $\CuspDt$ equals the discrete series subspace $L_d^2(X\col\tau)$
of $L^2(X\col\tau)$
and that $\TDt\colon\CciXt\to\CuspDt$ is the restriction of the orthogonal
projection $L^2(X\col\tau)\to L_d^2(X\col\tau)$.
In particular, the objects $\CuspDt$ and $\TDt$ are independent of the
choice of the residue weight $t$.
In the present paper $t$ is fixed throughout and
we do not need these properties. However, 
to simplify notation let
$\TD\coleq\TDt$ and $\CuspD\coleq\CuspDt.$
\end{rem} 

Fix $F\subset\Delta$ and let $\faqF\subset\faq$ be defined as
in [6, p.\ 41]. For each
$v\in\cW$ let $$\XFv=M_F/M_F\cap vHv^{-1}$$ 
be the reductive symmetric space defined as in [6, p.\ 51]. 
We use the  notation of [6, pp.\ 51, 52], related to this space.
Put $\tau_F=\tau|_{M_F\cap K}$ and let
the finite dimensional space  
$$
\CFv=\cA^{\stt}(\XFv\col\tau_F)
\subset\Cinf(\XFv\col\tau_F)$$ 
be the analog for $\XFv$ of the space $\CuspD$ of (5.2); cf.\ 
[6, eq.\ (10.7)], where the space is denoted $\cC_{F,v}$.
The assumption made before (5.2), that the vectorial part of $X$
vanishes,
holds for $\XFv$. For $\psi\in\CFv$ we have defined the generalized 
Eisenstein integral $\nE_{F,v}(\psi\col\nu)\in\Cinf(X\col\tau)$
in [6, Def.\ 10.7]; it is a linear function of $\psi$ and a
meromorphic function of $\nu\in\faqFdc$. Let us recall the definition.

The space $\CFv$ is spanned by elements $\psi\in\Cinf(\XFv\col\tau_F)$
of the form 
\begin{equation}
\psi(m)=\psi_{y,u}(m)=\stKF(\XFv\col m\col y)u
\end{equation}
for some $y\in\XFvp$, $u\in V_\tau$. Here 
$\stKF(\XFv\col\,\cdot\,\col\,\cdot\,)$ 
is the analog for $\XFv$ of the kernel $\KDt$, 
the residue weight $\stt\in\WT(\gS_F)$ is defined in [5, eq.\ (3.16)].
By definition
\begin{equation}
\nEFv(\psi_{y,u}\col\nu\col x)=\sum_{\gl\in\gL(\XFv,F)} 
\Res^{\stP,\stt}_\gl \bigl[
\nE(\nu-\,\cdot\,\col x)\after\inj_{F,v}
\Eastp(\XFv\col-\,\cdot\,\col y)u
\bigr]
\end{equation}
for $x\in X$.
Here $\Eastp(\XFv\col\gl\col y)=\Ep(\XFv\col-\bar\gl\col y)^*$ and  
$\gL(\XFv,F)\subset\faqdFp$ is the set defined
in [6, eq.\ (8.7)].
The generalized Eisenstein integral $\nEFv(\psi\col\nu\col x)$ 
is defined for $\psi\in\CFv$ by (5.4) and linearity; 
the fact that it is well defined
is shown in [6, Lemma 10.6], by using the induction of relations of
[7].
Let 
\begin{equation}
\psi=\sum_v\psi_v\in\CF\coleq\oplus_{v\in\cWF}\,\,\CFv,
\end{equation}
where $\cWF$ is as in [6, above Lemma 8.1]. Define
\begin{equation}
\nE_F(\psi\col\nu\col x)=\sum_{v\in\cWF}\nEFv(\psi_v\col\nu\col x).
\end{equation}

\begin{rem} {\it A priori\/} the generalized Eisenstein integral
$\nE_F(\psi\col\nu\col x)$
depends on the choice of the residue weight
$t$. In fact, already the parameter space $\CFv$ for $\psi$ depends 
on $t$ through the residue weight $\stt$. However, according to Remark 5.2{} 
(applied to the symmetric space $\XFv$) the latter 
is actually not the case. Once the independence of $\CFv$ on 
$\stt$ has
been established, it follows from the characterization in [8, Thm.\ 9.3],
that $\nE_F(\psi\col\nu\col x)$ is independent of $t$.
Therefore, this parameter is not indicated in the notation.
The independence of $t$ is not used in the present paper.
\end{rem}

\begin{lem} Let $\psi=\psi_{y,u}\in\CFv$ be given by\/
{\rm (5.3)} with $y\in\XFv${\rm ,} $u\in V_\tau$. Then
\begin{multline}
\nEFv(\psi_{y,u}\col\nu\col x)\\
=\sum_{\gl\in\gL(\XFv,F)} 
\Res^{\stP,\stt}_\gl \bigl[ \sum_{s\in W^F}
\Eps(\nu+\,\cdot\,\col x)\after\inj_{F,v}
\East(\XFv\col\,\cdot\,\col y)u
\bigr]
\end{multline}
for $x\in \Xp$ and generic $\nu\in\faqFdc$.
\end{lem}

\Proof If $y\in\XFvp$ then (5.4) holds and (5.7)
follows from [6, eq.\ (8.9]). The map $y\mapsto\psi_{y,u}$,
$\XFv\to\CFv$ is continuous, and $\nEFv(\psi\col\nu\col x)$ is linear
in $\psi$, hence the left side of (5.7) is continuous in $y\in\XFv$.
The other side is continuous as well, so (5.7) follows by the 
density of $\XFvp$ in $\XFv$.\Endproof\vskip4pt 

Let 
$$f\mapsto\TstF(\XFv\col f),\quad
\Cci(\XFv\col\tau)\to\CFv \subset\Cinf(\XFv\col\tau)$$ be the analog for $\XFv$
of the operator $\TD$ of (5.1) (with respect to some choice of
invariant measure $dy$ on $\XFv$). 
The operator
$\TstF(\XFv\col f)$ should not be confused with the operator $\TFt$ of 
[6, eq.\ (5.5)],
which maps between function spaces on $X$.
In the following lemma we examine the generalized Eisenstein
integral $\nE_{F,v}(\TstF(\XFv\col f)\col\nu)$. Let 
the Fourier transform associated with $\XFv$ be denoted
$f\mapsto\Fou(\XFv\col f)$. It maps 
$\Cci(\XFv\col\tau)$ into $\Mer(\faqFdperpc,\gS_F)\otimes\oCFv$ and is
given by (see (2.1))
\begin{equation}
\Fou(\XFv\col f)(\nu)=\int_{\XFv} \East(\XFv\col\nu\col y)f(y)\,dy,
\quad(\nu\in\faqFdperpc).
\end{equation}

\begin{lem}
Let $f\in\Cci(\XFv\col\tau)$ and let
$\psi=|W_F|^{-1}\TstF(\XFv\col f)\in\CFv$.
Then 
\begin{equation}
\nE_{F,v}(\psi\col\nu\col x)=
\sum_{\gl\in\gL(\XFv,F)} \Res^{\stP,\stt}_\gl \bigl[\sum_{s\in W^F}
\Eps(\nu+\,\cdot\,\col x)\after\inj_{F,v}
\Fou(X_{F,v}\col f)
(\,\cdot\,)\bigr]
\end{equation}
for $x\in\Xp$ and generic $\nu\in\faqFdc$.
\end{lem}

\Proof 
For each $y\in\XFv$ let $\psi_y\in\Cinf(\XFv\col\tau)$ be defined
by
$\psi_y(m)=\psi_{y,f(y)}(m)=\stKF(\XFv\col m\col y)f(y)
$;
cf.\ (5.3).
Then $\psi_y\in\CFv$ and $y\mapsto\psi_y$ is continuous
into this space. We conclude from (5.1), applied to $\XFv$, 
that $\psi=\int_{\XFv}\psi_y\,dy$ pointwise on $\XFv$, and hence also
as a $\CFv$-valued integral.
The Eisenstein integral $\nEFv(\psi\col\nu\col x)$ 
is linear in the first variable, hence
we further conclude that
\begin{equation}
\nEFv(\psi\col\nu\col x)=\int_{\XFv}\nEFv(\psi_y\col\nu\col x)\,dy.
\end{equation}
It follows from Lemma 5.4{} that
\begin{multline*}
\nEFv(\psi_y\col\nu\col x)\\
=\sum_{\gl\in\gL(\XFv,F)} 
\Res^{\stP,\stt}_\gl \bigl[ \sum_{s\in W^F}
\Eps(\nu+\,\cdot\,\col x)\after\inj_{F,v}
\East(\XFv\col\,\cdot\,\col y)f(y)
\bigr]
\end{multline*}
for $x\in\Xp$.
We insert this relation into (5.10) and take
the residue
operator outside
the integral over $y\in\supp f\subset\XFv$. The justification is
similar to that given in the proof of Lemma 3.8.
Using (5.8) we then obtain (5.9).\hfq

\begin{lem} The expressions {\rm (5.4), (5.7), (5.9)} 
remain valid if the set of summation
$\gL(\XFv,F)$ is replaced by any finite subset $\gL$ of $\faqFdperp$
containing $\gL(\XFv,F)$.\end{lem}

\Proof It follows from [6, Lemma 10.6], that the sum in (5.4) 
remains unchanged if $\gL(\XFv,F)$ is replaced by $\gL$.
That the same conclusion holds for (5.7) and (5.9) is then seen
as in the proofs of Lemmas 5.4{} and 5.5.
\hfq

\section{Induction of Arthur-Campoli relations}

In this section we prove in Theorem 6.2{} a result that will play a 
crucial role for the Paley-Wiener theorem. It shows that Arthur-Campoli 
functionals on the smaller symmetric space $\XFv$ induce  
Arthur-Campoli functionals on the full space $X$.
The result is established by means of the theory of 
induction of relations developed in [7, Cor.\ 16.4].
The corresponding result in the group case is [1, Lemma III.2.3],
however, for the unnormalized Eisenstein integrals.
Let $F\subset\Delta$,
and let $S\subset\faqFdperpc$ be finite.  

\begin{lem} \hskip-5pt Let 
$\Hyp$ be a  
$\gS$-configuration in $\faqdc${\rm ,} and let 
$\cL\!\in\!\Mer(\faqFdperpc,\gS_F)^*_\laur$ with $\supp\cL\subset S$.
\begin{itemize}
\ritem{\rm (i)} The set of affine hyperplanes in $\faqcdF${\rm ,}
$$\Hyp_F(S)=\cup_{a\in S}\,\{H'\mid\exists H\in\Hyp\colon\
a+H'=(a+\faqcdF)\cap H\subsetneq a+\faqcdF\},$$
is a $\gS_r(F)$-configuration{\rm ,}
which is real if $\Hyp$ is real and $S\subset\faqFdperp$.
The corresponding set of regular points is
$$\hskip-14pt\reg(\faqcdF,\Hyp_F(S))=\{\nu\in\faqcdF\mid 
\forall a\in S, H\in\Hyp: a+\nu\in H\Rightarrow a+\faqcdF\subset H\}.$$

\ritem{\rm (ii)} For each $\gf\in\Mer(\faqdc,\Hyp)$ and each 
$\nu\in\reg(\faqFdc,\Hyp_F(S))$
there exists a neighborhood $\Omega$ of $S$ in $\faqFdperpc$ such that
the function $\gf^\nu\colon\gl\mapsto\gf(\gl+\nu)$ belongs to 
$\Mer(\Omega,\gS_F)$.

\ritem{\rm (iii)} 
Fix $\nu\in\reg(\faqFdc,\Hyp_F(S))$.
There exists a
Laurent functional \/{\rm (}\/in general not unique\/{\rm )}\/
$\cL'\in\Mer(\faqdc,\gS)^*_\laur${\rm ,}
supported by the set $\nu+S${\rm ,} such that $\cL'\gf=\cL\gf^\nu$
for all $\gf\in\Mer(\faqdc,\Hyp)$.

\ritem{\rm (iv)} The function $\cL_*\gf\colon\nu\mapsto\cL\gf^\nu$ 
belongs to $\Mer(\faqFdc,\Hyp_F(S))$ for each $\gf\in\Mer(\faqdc,\Hyp)$.

\ritem{\rm (v)} The map $\cL_*$
maps $\Mer(\faqdc,\Hyp)$ continuously
into $\Mer(\faqFdc,\Hyp_F(S))$
and if $\Hyp$ is real{\rm ,}  
$\Pal(\faqd,\Hyp)$ continuously
into $\Pal(\faqFd,\Hyp_F(S))$.
\end{itemize}
\end{lem}

\Proof See [7, Cor.\ 11.6 and Lemma 11.7]. The continuity
in (v) between the $\Mer$ spaces is proved in [7, Cor.\ 11.6(b)];
the continuity between the $\Pal$ spaces is similar, see also
[5, Lemma 1.10].
\Endproof\vskip4pt 

Let $\Hyp=\Hyp(X,\tau)$ and let
$\nu\in\reg(\faqFdc,\cH_F(S))$.
Let $v\in\cWF$ and let $\pr_{F,v}\colon
\oC\to\oC_{F,v}$ be the projection operator defined by
[7, (15.3)]. 

\begin{theorem}
For each $\cL\in\ACXFv$ with $\supp\cL\subset S$ 
there exists a Laurent functional \/{\rm (}\/in general not unique\/{\rm )}\/
$\cL'\in\ACX${\rm ,}
supported by the set $\nu+S${\rm ,} such that
\begin{equation}
\cL[\pr_{F,v}\gf(\nu+\,\cdot\,)]=\cL'\gf,
\end{equation}
for all $\gf\in\Mer(\faqd,\Hyp)\otimes\oC$. 
In particular{\rm ,} if in addition $S\subset\faqFdperp$
then 
\begin{equation}
\cL[\pr_{F,v}\gf(\nu+\,\cdot\,)]=0
\end{equation}
for all $\gf\in\PACX$.
\end{theorem}

\Proof  The existence of $\cL'\in\LauoC$ such that
(6.1) holds follows from Lemma 6.1 (iii). We will show
that every such element $\cL'$ belongs to $\ACX$.
If $\nu\in\reg(\faqFd,\cH_F(S))$
the statement (6.2) is then straightforward from the
definition of $\PACX$, and in general
it follows by meromorphic continuation.

That $\cL\in\ACXFv$ means by definition that
it belongs to $$\Mer(\faqFdperpc,\gS_F)^*_\laur\otimes\oC_{F,v}^*$$ 
and satisfies
\begin{equation}
\cL[\East(\XFv\col\,\cdot\,\col m)u]=0
\end{equation}
for every $m\in\XFv$, $u\in\Vtau$. By (6.1) the claim that
$\cL'\in\ACX$ amounts to
\begin{equation}
\cL[\pr_{F,v}\East(X\col\nu+\,\cdot\,\col x)u]=0
\end{equation}
for all $x\in X$.
This claim will now be established by means of [7, Cor.\ 16.4].\vskip2pt

If $\psi \in \Mer(\faqFdperpc,\gS_F),$ then the function 
$\psi^\vee\colon \gl \mapsto
\overline{\psi(-\bar\gl)}$ belongs to\break $\Mer(\faqFdperpc,\gS_F)$ as well.
If $\cL \in \Mer(\faqFdperpc,\gS_F)^*_\laur,$ then it is readily seen that 
there exists a unique $\cL^\vee \in \Mer(\faqFdperpc,\gS_F)^*_\laur$ such that
\begin{equation}
\Lau^\vee \psi = (\Lau \psi^\vee)^*
\end{equation}
for all $\psi \in \Mer(\faqFdperpc,\gS_F);$ here the superscript 
$*$ indicates that the complex conjugate is taken.
The maps $\psi\mapsto\psi^\vee$ and  $\Lau\to\Lau^\vee$ are antilinear.
More generally, if $\Hil$ is a Hilbert space and $v \in \Hil,$ 
then by $v^*$ we denote the element
of the dual Hilbert space $\Hil^*$ defined by $v^*: w \mapsto \inp{w}{v}.$ 
The maps $(\psi,v)\mapsto \Psi^\vee\otimes v^*$ and
$(\Lau, v) \mapsto \Lau^\vee \otimes v^*$ induce antilinear 
maps from $\Mer(\faqFdperpc,\gS_F)\otimes \Hil $ 
to $\Mer(\faqFdperpc,\gS_F) \otimes \Hil^*,$ and from
$\Mer(\faqFdperpc,\gS_F)^*_\laur \otimes \Hil $ 
to $\Mer(\faqFdperpc,\gS_F)^*_\laur \otimes \Hil^*,$ 
which we denote by $\psi\mapsto\psi^\vee$ and
$\Lau \mapsto \Lau^\vee$ as well.
With this notation formula 
(6.5) is valid for all 
$\psi \in \Mer(\faqFdperpc,\gS_F) \otimes \Hil \otimes \Vtau$
and all
$\Lau \in \Mer(\faqFdperpc,\gS_F)^*_\laur \otimes \Hil.$ 
It is then an identity between members of $\Vtau$.

Notice that by definition of $\East(\XFv\col\,\cdot\,\col m)$
it is the $\psi^\vee$ of $$\psi=\nE(\XFv\col\,\cdot\,\col m)\in
\Mer(\faqFdperpc,\gS_F) \otimes\oC_{F,v}^*\otimes \Vtau.
$$
It now follows from (6.5) and (6.3) that 
\begin{equation}
\Lau^\vee(\nE(\XFv\col\,\cdot\,\col m))= 0
\end{equation}
for all $m \in \XFv,$ with $\Lau^\vee\in\Mer(\faqFdperpc,\gS_F)^*_\laur
\otimes\oC_{F,v}$
defined as above. 
Let
$$\Lau_2=(1\otimes\inj_{F,v})\Lau^\vee
\in \Mer(\faqFdperpc,\gS_F)^*_\laur
\otimes\oC,$$
then $
\Lau_2(\nE(\XFu\col\,\cdot\,\col m)\after\pr_{F,u})= 0
$ for all $u\in\cWF$, by (6.6) and [7, (16.2)].
In view of [7, Cor.\ 16.4] with $\Lau_1=0$ this implies that
\begin{equation}
\Lau_2[\nE(X\col\nu+\,\cdot\,\col x)]= 0
\end{equation}
for $x\in X_+$, hence by continuity also for $x\in X$.
Since 
$\cL_2=(\cL(1\otimes\pr_{F,v}))^\vee$
we readily obtain (6.4) by
application of (6.5) to (6.7). 
\hfq

\section{A property of the Arthur-Campoli relations}

The aim of this section is to establish a result, Lemma 
7.4, which elaborates on the definition of the space
$\ACX$ by means of some
simple linear algebra.

For any finite set $S\subset\faqdc$ we denote by $\cOS$ the space of germs 
at $S$ of functions $\phi\in\cO(\Omega)$, holomorphic on 
some open neighborhood $\Omega$ of $S$. Moreover, if
$\Omega$ is an open neighborhood of $S$ and
$d\colon\gS\to\Natural$ a map, then by
$\Mer(\Omega,S,\gS,d)$ we denote the space of meromorphic functions
$\psi$ on $\Omega$, whose germ at $a$ belongs to
$\pi_{a,d}^{-1}\cO_a$ for each $a\in S$. Here
$$\pi_{a,d}(\gl)=\Pi_{\alpha\in\Sigma}\,\langle\alpha,\gl-a\rangle^{d(\alpha)}$$
for $\gl\in \faqdc$ (cf.\ [7, eq.\ (10.1)]). Finally, we put
$\Mer(\Omega,S,\Sigma)=\cup_d \Mer(\Omega,S,\gS,d).$

\begin{lem}
 Let\/ $\lsp\subset\LauoC$ be a finite dimensional linear 
subspace{\rm ,} and let $S$ denote the finite set\/
$\supp\lsp\coleq\cup_{\cL\in\lsp}\supp\cL\subset\faqcd$. 
Then there exists a finite dimensional linear subspace $V\subset\CciXt$
with the following properties\/{\rm :}\/
\begin{itemize} 

\item[\rm (i)] Let $\Omega\subset\faqcd$ be an open neighborhood of $S$ 
and let $\psi\in\Mer(\Omega,S,\Sigma)\otimes\oC$ be annihilated 
by $\lsp\cap\ACX$. Then there exists a unique
function $f=f_\psi\in V$ such that $\cL\Fou f=\cL\psi$ for all $\cL\in\lsp$.
\vglue-32pt
\phantom{up}
\item[\rm (ii)] The map $\psi\mapsto f_\psi$ has the following form. 
There exists a $\Hom(\oC,V)$-valued Laurent functional 
$\cL'\in\lsp\otimes V\subset\Laufu\otimes\Hom(\oC,V)$ such that 
$f_\psi=\cL'\psi$  for all~$\psi$. 
\end{itemize}
\end{lem}
 
\phantom{up} \vglue-22pt
We first formulate a result in linear algebra, and then deduce the above
result.
 \vglue-20pt \phantom{up}
\begin{lem} Let $A${\rm ,} $B$ and $C$ be linear spaces with $\dim C<\infty${\rm ,}
and let $\alpha\in\Hom(A,B)$ and $\beta\in\Hom(B,C)$ be given. 
Put $C'=\beta(\alpha(A))$.
Then there exists a finite dimensional linear subspace $V\subset A$
with the property that{\rm ,} for each $\psi\in\beta^{-1}(C')${\rm ,} 
there exists a unique
element $f_\psi\in V$ such that $\beta(\alpha(f_\psi))=\beta(\psi)$.
Moreover{\rm ,} there exists an element $\mu\in\Hom(C,V)$ such that 
$f_\psi=\mu(\beta(\psi))$ for all $\psi$.
\end{lem}

\phantom{up}
\vglue-20pt

{\it Proof}.  The proof is shorter than the statement.
Since $\beta\circ\alpha$ maps $A$ onto $C'$ we can
choose $V\subset A$ such that the restriction of 
$\beta\circ\alpha$ to it is bijective $V\to C'$. Then $f_\psi\in V$
is uniquely determined by $\beta\circ\alpha(f_\psi)=\beta(\psi)$,
and if $\mu\colon C\to V$ is any linear extension of
$(\beta\circ\alpha)^{-1}\colon C'\to V$, the relation
$f_\psi=\mu(\beta(\psi))$ holds for all $\psi$.
\hfq

\vskip8pt {\it Proof of Lemma {\rm 7.1}}.  It is easily seen by using a basis for $\lsp$ that $S$ is a 
finite set.

We shall apply Lemma 7.2{} with $A=\CciXt$, 
$B=\Mer(\Omega,S,\Sigma)\otimes\oC$
and $C=\lsp^*$, the linear dual of $\lsp$. 
Furthermore, as $\alpha\colon A\to B$ we use the Fourier transform $\Fou$
followed by taking restrictions to $\Omega$, 
and as $\beta\colon B\to C=\lsp^*$
we use the map induced by the pairing $(\cL,\psi)\mapsto\cL\psi$,
$\cL\in\lsp$, $\psi\in B$. 

We now determine the image $C'=\beta(\alpha(A))$. By definition it consists
of all the linear forms on $\lsp$ given by the application of 
$\cL\in\lsp$ to a function in $\Fou(\CciXt)$. Hence the polar subset
$C'{}^\perp\subset\lsp$ is exactly the set of $\cL\in\lsp$ that 
annihilate  $\Fou(\CciXt)$. By Lemma 3.8, an element 
$\cL\in\lsp$ annihilates $\Fou(\CciXt)$ if and only if it belongs to 
$\ACX$. Hence $C'{}^\perp=\lsp\cap\ACX$. Thus $\beta^{-1}(C')$
consists precisely of those elements 
$\psi\in B=\Mer(\Omega,S,\Sigma)\otimes\oC$ that are annihilated
by $\lsp\cap\ACX$.

The lemma now follows immediately from Lemma 7.2.
\hfq 
\vglue-22pt
\phantom{up}
\begin{lem}
 Let $\cL\in\Laufu$  and let $\phi\in\cOS$ where $S=\supp\cL$.
The map $\cL_\phi\colon\psi\mapsto\cL(\phi\psi)$ is a Laurent 
functional in $\Laufu${\rm ,} supported at $S$. \end{lem}

\phantom{up}
\vglue-22pt
{\it Proof}.
(See also [7, eq.\ (10.7)].)
For each $a\in S$, let $u_a=(u_{a,d})$ be the string that represents 
$\cL$ at $a$. 
Let $\Omega$ be an open neighborhood
of $S$. Fix $d\colon\gS\to\Natural$.
For 
$\psi\in\Mer(\Omega,S,\gS,d)$ we have 
$\cL_\phi\psi=\sum_{a\in S} u_{a,d}[\pi_{a,d}\phi\psi](a).$
Hence by the Leibniz rule we can write
\begin{equation}
\cL_\phi\psi=\sum_{a\in S}
\sum_i u^1_{a,i}[\phi](a)\, u^2_{a,i}[\pi_{a,d}\psi](a)
\end{equation}
for finitely many $u_{a,i}^1, u_{a,i}^2\in S(\faqd)$. Thus $\cL_\phi$ has 
the form 
required of a Laurent functional with support in $S$.
\hfq

\begin{lem} Let $\cL_0\in\Laufu$ and let $d\colon\gS\to\Natural$. 
There exists a finite dimensional linear subspace 
$V\subset\CciXt$ with the following  properties\/{\rm : }\/ 
\begin{itemize}
\ritem{\rm (i)} Let $\Omega\subset\faqcd$ be an open neighborhood of 
$S\coleq\supp\cL_0$ 
and let $\psi\in\Mer(\Omega,S,\gS,d)\otimes\oC$. Assume that
$\cL\psi=0$ for all $\cL\in\ACX$ with $\supp\cL\subset S$.
Then there exists a unique function $f=f_\psi\in V$ such that 
$\cL_0(\phi\Fou f)=\cL_0(\phi\psi)$ for all $\phi\in\cOS\otimes\oC^*$. 
\vglue-18pt
\phantom{up}
\ritem{\rm (ii)}  The map $\psi\mapsto f_\psi$ has the following form. 
There exists 
a $Hom(\oC,V)$-valued germ $\phi'\in\cOS\otimes\Hom(\oC,V)$ such that 
$f_\psi=\cL_0(\phi'\psi)$  for all $\psi$.\end{itemize}
\end{lem}

\Proof 
We may assume that the given $d\in\Natural^{\gS}$ satisfies the requirement  that 
$\Fou f|_\Omega$ 
belongs to $\Mer(\Omega,\gS,d)\otimes\oC$ for all $f\in\CciXt$, for some
neighborhood $\Omega$ of $S$ (otherwise we 
just replace $d$ by a suitable successor in $\Natural^{\gS}$).

Let $\cOSa=\cOS\otimes\oC^*$ and
let $\cOSo$ denote the subspace of $\cOSa$ consisting of the elements 
$\phi\in\cOSa$ for which the Laurent functional $\cL_{0\phi}\colon
\psi\mapsto\cL_0(\phi\psi)$ in $\Laufu\otimes\oC^*$
annihilates $\Mer(\Omega,S,\gS,d)\otimes\oC$ (with the fixed element $d$), 
for all neighborhoods $\Omega$ of $S$. It follows immediately from (7.1), 
applied componentwise on $\oC$, that an element $\phi\in\cOSa$ belongs to
$\cOSo$ if a finite number of fixed linear forms on $\cOSa$ annihilate it; 
hence $\dim \cOSa/\cOSo<\infty$. 
Fix a complementary subspace $\cOSp$ of $\cOSo$ in $\cOSa$, and let 
$$\lsp=\{\cL_{0\phi}\mid\phi\in \cOSp\}\subset\LauoC.$$ 
Choose $V\!\subset\!\Cci(X\!\col\!\tau)$ according to Lemma 7.1. Then for 
each $\psi\!\in\!\Mer(\Omega,S,\gS,d)\break\otimes\oC$ satisfying $\cL\psi=0$
for all $\cL\in\lsp\cap\ACX$, there exists a unique function
$f_\psi\in V$ such that $\cL\Fou f_\psi=\cL\psi$ for all $\cL\in\lsp$.
Thus $\cL_0(\phi\Fou f_\psi)=\cL_0(\phi\psi)$ for all 
$\phi\in\cOSp$, and this property determines $f_\psi$ uniquely.
On the other hand, by the definition of $\cOSo$ we have
$\cL_0(\phi\Fou f_\psi)=0=\cL_0(\phi\psi)$
for
$\phi\in\cOSo$. 
Thus
$\cL_0(\phi\Fou f_\psi)=\cL_0(\phi\psi)$ holds for all $\phi\in\cOSa$.

The statement (ii) follows immediately 
from the above and the corresponding statement in Lemma 7.1.
\hfq

\section{Proof of Theorem 4.4}

The inversion formula for the Fourier transform that was obtained
in [6, Thm.\ 1.2], reads
\begin{equation}
f(x)=\cT\Fou f(x)=\sum_{F\subset\Delta}\TFt f(x),\qquad
x\in\Xp,
\end{equation}
where the term in the middle is the pseudo wave packet (4.1)
and where the operators on the right-hand side
are as defined in [6, eq.\ (5.5)]. Motivated by
the latter definition we define,
for $F\subset\Delta$, $\gf\in\cP(X\col\tau)$ and $x\in\Xp$, 
\begin{multline}
\cTF\gf(x)
=|W|\, t(\faqFp)\\
\cdot
\int_{\geps_F+i\ayFd}\sum_{\gl\in\gL(F)} \Res^{P,t}_{\gl+\faqdF} 
\bigl[\sum_{s\in W^F} \Eps(\,\cdot\,\col x)
\gf(\,\cdot\,)\bigr](\gl+\nu)\,d\mu_{\faqdF}(\nu)
\end{multline}
so that
$\TFt f=\cTF\Fou f$. The element
$\geps_F\in\fa^{*+}_{F\iq}$, the set $\gL(F)\subset\faqdFp$
and the measure
$d\mu_{\faqdF}$ on $i\ayFd$ are as defined in [6, p.\ 42]
(with $\Hyp$ equal to the union of $\Hyp(X,\tau)$ with the set
of singular hyperplanes for
$\Ep$).
It follows from\break [6, eq.\ (4.2)]
and [5, Lemma 1.11], that the integral in (8.2) converges, and
that $\cTF\gf\in\Cinf(\Xp \col\tau)$. Moreover,
\begin{equation}
\cT\gf=\sum_{F\subset\Delta}\cTF\gf,
\end{equation}
in analogy with the second equality in (8.1);
see the arguments leading up to [6, eq.\ (5.3)].

The existence of a smooth extension of $\psW\gf$ will be proved 
by showing that $\cTF\gf$ has the same property, for each~$F$. 
We shall do this by exhibiting it
as a wave packet of generalized Eisenstein integrals. 

Let $\Hyp$ denote the union of $\Hyp(X,\tau)$ with
the set of all affine hyperplanes in $\faqdc$
along which $\gl\mapsto\Eps(\gl\col x)$ is singular, 
for some $x\in\Xp$, $s\in W$.
By Lemma 2.1 this is a real $\gS$-configuration
and there exists $d\maps\Hyp\to\Natural$ such that
$\Eps(\,\cdot\,\col x)\in\Mer(\faqd,\Hyp,d)\otimes\Hom(\oC,V_\tau)$
for all $x\in\Xp$ and $s\in W$. 

\begin{lem} Let $F\subset\Delta$ 
and $v\in\cWF$.  Let $\cL\in\Mer(\faqFdperpc,\gS_F)^*_\laur$
with $S\coleq\supp\cL\subset\faqFdperp$. 
There exist a finite dimensional 
linear subspace $V\subset\Cci(\XFv\col\tau)$ 
and for each $\nu\in\reg(\faqFdc,\cH_F(S))$ a linear 
map $\gf\mapsto f_{\nu,\gf}${\rm ,} $\PACX\to V$, such that 
\begin{multline}
 \cL\Bigl[\sum_{s\in W^F} \Eps(\nu+\,\cdot\,\col x)\after\inj_{F,v}\after
\pr_{F,v}
\gf(\nu+\,\cdot\,)\Bigr]
\\
 =\cL\Bigl[\sum_{s\in W^F} \Eps(\nu+\,\cdot\,\col x)\after\inj_{F,v}
\Fou(X_{F,v}\col f_{\nu,\gf})(\,\cdot\,)\Bigr]
\end{multline}
for all $x\in\Xp$.

Moreover{\rm ,} the elements $f_{\nu,\gf}\in V$ can be chosen of the 
following form. 
There exists a Laurent functional
$\cLa\in\Mer(\faqFdperpc,\gS_F)^*_\laur\otimes\Hom(\oC_{F,v},V)${\rm ,} 
supported by $S${\rm ,} such that
\begin{equation}
 f_{\nu,\gf}=\cLa[\pr_{F,v}\gf(\nu+\,\cdot\,)]
\end{equation}
for all $\nu\in\reg(\faqFdc,\cH_F(S))$  and all $\gf\in\PACX$.
\end{lem}

\Proof 
For each $\nu\in\reg(\faqFdc,\cH_F(S))$ and $a\in S$ the element
$a+\nu$ is only contained in a given hyperplane from $\Hyp$ if this hyperplane 
contains all of $a+\faqFdc$. Let $\Hyp(a+\faqFdc)$ denote the (finite) 
set of such 
hyperplanes, and let $\Hyp(S+\faqFdc)=\cup_{a\in S} \Hyp(a+\faqFdc)$.
Let $d\colon\Hyp\to\Natural$ be as mentioned before the lemma,
and let the polynomial function $p$ be given by (2.6)
with $\omega=\nu+S$, where $\nu\in\reg(\faqFdc,\cH_F(S))$. Then
$$p(\gl)=\prod_{H\in\Hyp(S+\faqFdc)}
(\langle\ga_H,\gl\rangle-s_H)^{d(H)},$$
and thus $p$ is independent of $\nu$.
Moreover, since $a+\faqFdc\subset H$ we conclude that $\ga_H\in\gS_F$
for all $H\in\Hyp(S+\faqFdc)$. Hence $p(\nu+\gl)=p(\gl)$ for
$\nu\in\faqFdc$ and $\gl\in\faqFdperpc$.
The maps
$$\gl\mapsto p(\gl)\Eps(\nu+\gl\col x),\quad
\faqFdperpc\to\Hom(\oC,V_\tau),$$
are then holomorphic at $S$ for all $\nu\in\reg(\faqFdc,\cH_F(S))$, 
$s\in W$ and 
$x\in\Xp$. 

Choose $d_0\in\Natural$ such that $d_{X,\tau}(H)\le d_0$ for all
$H\in\Hyp(S+\faqFdc)\cap\Hyp(X,\tau)$ 
and define
$d'\maps\gS_F\to\Natural$ by $d'(\ga)=d_0$
for all $\ga$.
Then, for each $\nu\in\reg(\faqFdc,\cH_F(S))$ and 
$\gf\in\Mer(\faqd,\Hyp(X,\tau),d_{X,\tau})\otimes\oC$ 
the function 
$$\psi^{\nu,\gf}\coleq\pr_{F,v}\circ\varphi^\nu\colon
\gl\mapsto \pr_{F,v}\gf(\nu+\gl)$$
on $\faqFdperpc$ belongs to $\Mer(\Omega,\gS_F,d')\otimes\oC_{F,v}$
for some neighborhood $\Omega$ of $S$ (cf.\ Lemma 6.1). 
If in addition $\gf\in\PACX$ then
by Theorem 6.2{} this function is annihilated by
all elements of $\AC(\XFv\col\tau)$ supported by $S$.

Let $\cL_0$ be the functional on $\Mer(\faqFdperpc,\gS_F)$ defined by
$\cL_0\psi=\cL(p^{-1}\psi)$; it is easily seen that $\cL_0\in
\Mer(\faqFdperpc,\gS_F)^*_\laur$ and that $\supp\cL_0\subset S$.
Choose $V\subset\Cci(\XFv\col\tau)$ according to Lemma 7.4,
applied to $\XFv$, $\cL_0$ and $d'$.
Then there exists
for each $\nu\in\reg(\faqFdc,\cH_F(S))$  and  
$\gf\in\PACX$ a unique element 
$f_{\nu,\gf}=f_{\psi^{\nu,\gf}}\in V$ such that
$$\cL_0(\phi\Fou(\XFv\col f_{\nu,\gf}))=\cL_0(\phi\psi^{\nu,\gf})$$ 
for all 
$\phi\in
\cOS\otimes\oC^*_{F,v}$.
We apply this identity with 
$$\phi(\gl)=
p(\gl)\sum_{s\in W^F} \upsilon^*\after\Eps(\nu+\gl\col x)\after\inj_{F,v}$$ 
for arbitrary $\upsilon^*\in V_\tau^*$, and deduce (8.4).

According to Lemma 7.4{} (ii)
there exists $\phi'\in\cOS\otimes\Hom(\oCFv, V)$ such
that $f_{\nu,\gf}=\cL_0(\phi'\psi^{\nu,\gf})$.
The map $\cLa\colon\psi\mapsto\cL_0(\phi'\psi)$ 
is a $\Hom(\oCFv,V)$-valued
Laurent functional (see Lemma 7.3)
satisfying (8.5).  The linearity of $\gf\mapsto f_{\nu,\gf}$
follows from (8.5).
\hfq

\begin{lem}
Let $v\in\cWF$.
There exists a Laurent functional $$\cLFv\in 
\Mer(\faqFdperpc,\gS_F)^*_\laur\otimes\Hom(\oCFv,\CFv),$$
supported by the set $\gL\coleq\gL(F)\cup\gL(\XFv,F)${\rm ,} such that
\begin{multline}
\sum_{\gl\in\gL(F)} \Res^{P,t}_{\gl+\faqdF} \bigl[\sum_{s\in W^F}
\Eps(\,\cdot\,\col x)\after \inj_{F,v}\after\pr_{F,v}
\gf(\,\cdot\,)\bigr](\nu+\gl) \\
= 
\nE_{F,v}(\cLFv[\pr_{F,v}\gf(\nu+\,\cdot\,)]\col\nu\col x) 
\end{multline}
for all $\gf\in\PACX${\rm ,}
$x\in\Xp$ and generic $\nu\in\faqFdc$. Here{\rm ,}
generic means that
$\nu\in\reg(\faqFdc,\Hyp_F(\gL))${\rm ,} where $\Hyp$ is as defined
above Lemma {\rm 8.1}.
\end{lem}

\Proof 
In the expression on the left side of (8.6) 
we can replace the set $\gL(F)$ by $\gL$ 
(see [6, Lemma 7.5]).
Moreover, we can replace the residue operator 
$\Res^{P,t}_{\gl+\faqdF}$
by $\Res^{\stP,\stt}_\gl$ (see [6, eq.\ (8.5)]), 
which, as observed in [6, above eq.~(8.5)], can be
regarded as an
element in $\Mer(\faqFdperp,\gS_F)^*_\laur$, supported at $\gl$. We thus
obtain on the left of (8.6):
\begin{equation}
\sum_{\gl\in\gL} \Res^{\stP,\stt}_\gl \bigl[\sum_{s\in W^F}
\Eps(\nu+\,\cdot\,\col x)\after\inj_{F,v}\after\pr_{F,v} 
\gf(\nu+\,\cdot\,)\bigr].
\end{equation}
We obtain from Lemma 8.1  that there exist
a finite dimensional space $V\subset\Cci(\XFv\col\tau)$ 
and a Laurent functional $\cLa\in\Mer(\faqFdperp,\gS_F)^*_\laur\otimes
\Hom(\oCFv,V)$ 
supported by $\gL$, such that (8.7) equals 
\begin{equation}
\sum_{\gl\in\gL} \Res^{\stP,\stt}_\gl \bigl[\sum_{s\in W^F}
\Eps(\nu+\,\cdot\,\col x)\after \inj_{F,v}
\Fou(\XFv\col f_{\nu,\gf})(\,\cdot\,)\bigr].
\end{equation}
Here $f_{\nu,\gf}=\cLa[\pr_{F,v}\gf(\nu+\,\cdot\,)]\in V$
for $\nu\in\reg(\faqFdc,\Hyp_F(\gL))$.
We apply  Lemmas 5.5, 5.6  and obtain that (8.8) equals
$\nEFv(\psi\col\nu\col x)$ with $\psi=
|W_F|^{-1} \TstF(\XFv\col f_{\nu,\gf})\in\CFv$.

The map $f\mapsto|W_F|^{-1} \TstF(\XFv\col f)$ is linear 
$V\to\CFv$; composing it with  the coefficients of
$\cLa\in\Mer(\faqFdperp,\gS_F)^*_\laur\otimes\Hom(\oCFv,V)$ 
we  obtain  a Laurent functional $\cLFv\in
\Mer(\faqFdperp,\gS_F)^*_\laur\otimes\Hom(\oCFv,\CFv)$.  Now 
$\psi=\cLFv[\pr_{F,v}\gf(\nu+\,\cdot\,)]$,  and
(8.6) follows.
\hfq

\begin{theorem} Let $F\subset\Delta$.
There exists $$\cLb\in 
\Mer(\faqFdperpc,\gS_F)^*_\laur\otimes\Hom(\oC,\CF)$$
with support contained in $\gL(F)\cup[\cup_{v\in\cWF}\,\gL(\XFv,F)]${\rm ,} 
such that
\begin{equation}\cTF\gf(x)= \int_{\geps_F+i\ayFd}
\nE_F(\cLb[\gf(\nu+\,\cdot\,)]\col\nu\col x)
\,d\mu_{\faqFd}(\nu)
\end{equation}
for all $\gf\in\PACX${\rm ,} $x\in\Xp$. In particular{\rm ,} 
$\cTF\gf\in\Cinf(X\col\tau)${\rm ,} and $\gf\mapsto\cTF\gf$ is
continuous $\PACX\to\Cinf(X\col\tau)$.
\end{theorem}

\Proof
Recall, see (5.5) and [6, eq.\ (8.4)], that 
$$
\CF=\oplus_{v\in\cWF}\,\,\CFv,\quad
\oC=\oplus_{v\in\cWF}\, \inj_{F,v}\bigl(\oCFv\bigr). $$ 
Let $\cLFv$ 
be as in Lemma 8.2{} for each $v\in\cWF$,
and let $$\cLb=|W|t(\faqFp)\sum_{v\in\cWF} \cLFv\after\pr_{F,v}\in
\Mer(\faqFdperpc,\gS_F)^*_\laur\otimes\Hom(\oC,\CF).$$
The identity (8.9) then follows immediately from 
(8.2), (8.6), (5.6).
The remaining statements follow from Lemma 6.1(v)
combined with the estimate in [6, Lemma 10.8].\Endproof\vskip4pt 

As a corollary 
we immediately obtain (cf.\ (8.3)) that 
$\psW\gf\in C^\infty(X\col\tau)$ 
for  every $\gf\in\PACX$, and that
$\psW\colon\PACX\to\Cinf(X\col\tau)$ is continuous.
The proofs of Theorems 4.4{}, 
4.5{} and 3.6{} are then complete.

\section{A comparison of two estimates}

The purpose of this section is to compare the estimates (3.2) 
and (4.2), and to establish the facts mentioned in Remark 4.2.
The method is elementary Euclidean Fourier analysis.

Fix $R\in\real$ and let $\cQ=\cQ(R)$ denote the space of functions
$\phi\in\cO(\faqd(P,R))$ (see (2.7)) for which
\begin{equation}
\nu_{\omega,n}(\phi)\coleq
\sup_{\gl\in\omega+i\faqd}\,(1+ |\gl|)^{n} |\phi(\gl)|<\infty
\end{equation}
for all $n\in\Natural$ and all bounded sets
$\omega\subset\faqd(P,R)\cap\faqd$.
The space $\cQ$, endowed with the seminorms $\nu_{\omega,n}$,
is a Fr\'echet space. 

For $M>0$ we denote by $\cQ_M=\cQ_M(R)$ the subspace of $\cQ$ consisting of
the functions $\phi\in\cQ$ that satisfy the following:
For every strictly antidominant $\eta\in\faqd$ there exist
constants $t_\eta, C_\eta>0$ such that
\begin{equation}|\phi(\gl)|\le C_\eta (1+|
\gl|)^{-\dim\faq-1} e^{M|\Re\gl|}
\end{equation}
for all $t\geq t_\eta$ and $\gl\in t\eta+i\faqd$
(note that $t\eta+i\faqd\subset\faqd(P,R)$ for $t$ sufficiently large).

\def\glmu{\mu}
\begin{lem} {\rm (i)} Let $\gl_0\in\faqd(P,R)\cap \faqd$ and
let $\omega\subset\faqd(P,R)\cap \faqd$
be a compact neighborhood of $\gl_0$. Let $M>0$ and $N\in\Natural$.
There exist $n\in\Natural$ and $C>0$ such that
\begin{equation}
|\phi(\gl)|\le C(1+ |\gl|)^{-N} 
e^{M|\Re\gl|}\nu_{\omega,n}(\phi)
\end{equation}
for all $\gl\in\gl_0+\baraqd(P,0)$ and $\phi\in\cQ_M$.
\smallbreak 
 {\rm (ii)} $\cQ_M$ is closed in $\cQ$.
\smallbreak {\rm (iii)} Let $\phi\in\cQ_M$. Then $p\phi\in\cQ_M$ for
each polynomial $p$ on $\faqdc$.
\end{lem}

\Proof (i) {}From the estimates in (9.1) it follows
that $\glmu\mapsto\phi(\gl_0+\glmu)$ is a Schwartz function
on the Euclidean space $i\faqd$; in fact by a straightforward
application of Cauchy's integral formula we see that
every Schwartz-type seminorm of this function
can be estimated from above by (a constant times)
$\nu_{\omega,n}(\phi)$ for some $n$.
\par
Let $f\maps \faq\to\complex$ be defined by
\begin{equation}
f(x)=\int_{\gl_0+i\faqd} e^{\gl(x)} \phi(\gl) \,d\gl.
\end{equation}
Then $x\mapsto e^{-\gl_0(x)}f$ is a Schwartz function on $\faq$,
and by continuity of the Fourier transform
for the Schwartz topologies every Schwartz-seminorm of this function
can be estimated by one of the $\nu_{\omega,n}(\phi).$
Moreover, it follows from the Fourier inversion formula that
\begin{equation}
\phi(\gl)=\int_{\faq} e^{-\gl(x)}f(x) \,dx,
\end{equation}
for $\gl\in\gl_0+i\faqd$, where $dx$ is Lebesgue measure on $\faq$
(suitably normalized).
\par
It follows from (9.4) and
an application of Cauchy's theorem,
justified by (9.1), that $f(x)$ is independent of the choice of
the element $\gl_0$. Since this element was arbitrary in
$\faqd(P,R)\cap\faqd$,
we conclude that (9.5) holds for all
$\gl\in\faqd(P,R)$.

Let $\glmu\in\faqd(P,0)$ and let $\eta=\Re\glmu$. Then
$\eta$ is strictly antidominant. Let $t\ge t_\eta$.
Replacing $\gl_0$ by $t\eta$ in (9.4) and applying
(9.2) we obtain the estimate
$$
| f(x)|\leq  C_\eta \, e^{t\eta(x)}  e^{tM|\eta|} \int_{i\faqd}
(1+|\gl|)^{-\dim\faq-1} \, d\gl.
$$
By taking the limit as $t\to\infty$ we infer
that if $\eta(x)+M|\eta|<0$ then $f(x)=0$.
 
We use (9.5) to evaluate $\phi(\gl_0+\glmu)$. It follows from the
previous statement that we need only to integrate over the set where
$-\eta(x)\leq M|\eta|$. On
this set the integrand $e^{-(\gl_0+\glmu)(x)}f(x)$ is dominated by
$e^{M|\eta|}e^{-\gl_0(x)}|f(x)|$. Thus we obtain
\begin{equation}
|\phi(\gl_0+\glmu)| \le e^{M|\Re\glmu|}\int_{\faq} e^{-\gl_0(x)}|f(x)| \,dx
\end{equation}
for $\glmu\in\faqd(P,0)$,  hence, by continuity, also for
$\glmu\in\baraqd(P,0)$. Using (9.5) and partial integration, we obtain
a similar estimate   for
$\glmu(x_0)^k\phi(\gl_0+\glmu)$ for any $x_0\in\faq,$ $k\in\Natural$;
on the right-hand side of (9.6)
$e^{-\gl_0}f$ is then replaced by its $k$-th derivative in the
direction $x_0$. This shows that  for each $N\in\Natural$,
$(1+|\glmu|)^N|\phi(\gl_0+\glmu)|$
can be estimated in terms of $e^{M|\Re\glmu|}$ and a Schwartz-seminorm
of $e^{-\gl_0}f$. The latter seminorm may then be estimated by 
$\nu_{\omega,n}(\phi),$
for suitable $n$, and (9.3) follows, but with $\mu=\gl-\gl_0$
in place of $\gl$ on the right-hand side. Since 
$1+|\gl|\leq 1+|\gl_0|+|\mu|\leq (1+|\gl_0|)(1+|\mu|)$ and
$|\Re\mu|\leq|\Re\gl_0|+|\Re\gl|$, the stated form of (9.3) 
follows from that. 
\smallbreak
(ii) Let $\phi$ be in the closure of $\cQ_M$ in $\cQ$; then
by continuity (i) holds for $\phi$ as well.
Let $\eta$ be a given, strictly antidominant, element of $\faqd$.
Choose $t_\eta>0$ such that $\gl_0\coleq t_\eta\eta\in\faqd(P,R)$.
Now (9.2) follows from (9.3) with
$N=\dim\faq+1$.
Hence $\phi\in\cQ_M$.
\smallbreak
(iii) As before, let $\eta$ be given and choose $t_\eta>0$
such that $\gl_0=t_\eta\eta\in\faqd(P,R)$. Then by (i), (9.3)
holds, and since $N$ is arbitrary (9.2) follows with $\phi$
replaced by $p\phi$.
\hfq

\begin{lem} There exist a real $\gS$-configuration
$\Hypt${\rm ,} a map $\dt\colon\Hypt\to\Natural$ and a number $\geps>0$
with the following property.
Let $\gf\colon\faqdc\to\oC$ be any meromorphic function such that
\begin{itemize}
\ritem{\rm (i)} $\gf(s\gl)=\nC(s\col\gl)\gf(\gl)$ for all $s\in W$ 
and generic $\gl\in\faqdc${\rm ,}

\ritem{\rm (ii)} $\pi\gf$ is
holomorphic on a neighborhood of $\baraqd(P,0)$.
\end{itemize}
Then $\gf\in\Mer(\faqd,\Hypt,\dt)\otimes\oC$ and
$\pi\gf$ is holomorphic on
$\faqd(P,\geps)$. 
\end{lem}

Notice (cf.\ (2.3))
that (i), (ii) hold with $\gf=\East(\,\cdot\,\col x)v$,
for any $x\in X$, $v\in\Vtau$.
It follows that $\East(\,\cdot\,\col x)v\in
\Mer(\faqd,\Hypt,\dt)\otimes\oC$. Hence 
$\Hyp(X,\tau)\subset\Hypt$ and $d_{X,\tau}\preceq\dt|_{\Hyp(X,\tau)}$.

\Proof Let $\Hyp(X,\tau)$ and $d_{X,\tau}$
be as in Section 2, and
for each $s\in W$ let $\Hyp_s$, $d_s$ be such that
$\nC(s\col\,\cdot\,)\in\Mer(\faqd,\Hyp_s,d_s)$; cf.\ 
Lemma 2.1. Let 
$$\Hypt=\cup_{s\in W} \{sH\mid H\in\Hyp(X,\tau)\cup\Hyp_s\}.$$ Furthermore,
let $\dt\in\Natural^{\Hypt}$ be defined as follows.
We agree that $d_{X,\tau}(H)=0$ for $H\notin\Hyp(X,\tau)$ and
  $d_s(H)=0$ for $H\notin\Hyp_s$.
For $H\in\Hypt$ let 
$$\dt(H)=\max_{s\in W} d_{X,\tau}(s^{-1}H)+d_s(s^{-1}H).$$

We now assume that $\gf$ satisfies (i) and (ii).
Let $\gl_0\in\baraqd(P,0)$ and $s\in W$. Let $\pi_0$ denote the polynomial 
determined by (2.6) with $\omega=\{\gl_0\}$ and with
$\Hyp=\Hyp(X,\tau)$ and $d=d_{X,\tau}$. Since $\gl_0\in\baraqd(P,0)$,
we see that $\pi_0$ divides $\pi$ and the quotient $\pi/\pi_0$ is
nonzero at $\gl_0$. Hence $\pi_0\gf$ is holomorphic in a neighborhood
of $\gl_0$, by (ii). Likewise, let
$\pi_s$ denote the polynomial 
determined by (2.6) with $\omega=\{\gl_0\}$ and with
$\Hyp=\Hyp_s$ and $d=d_s$, then $\pi_s\nC(s\col\,\cdot\,)$
is holomorphic at $\gl_0$. Hence
$\pi_0\pi_s \nC(s\col\,\cdot\,)\gf$ is holomorphic at
$\gl_0$, and by (i) it follows that 
$\gl\mapsto\pi_0(s^{-1}\gl)\pi_s(s^{-1}\gl) \gf(\gl)$
is holomorphic at $s\gl_0$.
Let $\pi^\sim$ be defined by (2.6) 
with $\omega=\{s\gl_0\}$ and with $\Hyp=\Hypt$ and $d=\dt$.
Then the polynomial $\gl\mapsto\pi_0(s^{-1}\gl)\pi_s(s^{-1}\gl)$ divides
$\pi^\sim$, by the definition of $\dt$, and
hence $\pi^\sim\gf$ 
is holomorphic at $s\gl_0$. Since every point in $\faqdc$ can be written
in the form $s\gl_0$ with $\gl_0\in\baraqd(P,0)$ and $s\in W$, it
follows that $\gf\in\Mer(\faqd,\Hypt,\dt)\otimes\oC$.
The statement about the existence of $\geps$ is now an easy consequence
of (ii) and the local finiteness of $\Hypt$.
\Endproof\vskip4pt 

It follows from Lemma 9.2{} that a fixed number
$\geps$ can be chosen such that the condition in
(ii) of Definition 4.1 holds for all
$\gf\in\PXt$ simultaneously. In the following lemma, we fix such a number
$\geps>0$.

\begin{lem} Let $M>0$ and let $\omega\subset\faqd(P,\geps)$ 
be a compact neighborhood
of $0$. Let $N\in\Natural$. Then there exist $n\in\Natural$ and $C>0$
such that
\begin{equation}
\sup_{\gl\in\baraqd(P,0)} (1+|\gl|)^N e^{-M |\Re\gl|}
\|\pi(\gl)\gf(\gl)\| \le C\nu_{\omega, n}(\pi\gf)
\end{equation}
for all $\gf\in\PMXt$ \/{\rm (}\/see Definition {\rm 4.1)}. Moreover{\rm ,}
\begin{equation}
\PW_M(X\col\tau)=\Pal_M(X\col\tau)\cap\PACX,
\end{equation}
and this is a closed subspace of $\PACX$.\end{lem}

\Proof We first show that $\pi\gf\in \cQ_M(\geps)\otimes\oC$
for all $\gf\in\PMXt$. Let $\gf\in\PMXt$ and let
$R_1\in\real$ be sufficiently negative so that
$\gf$ is holomorphic on $\faqd(P,R_1)$.
Then $\gf\in\cQ_M(R_1)\otimes\oC$ and hence it follows
from Lemma 9.1 (iii) with $R=R_1$, 
applied componentwise to the $\oC$-valued function
$\gf$, that $\pi\gf\in\cQ_M(R_1)\otimes\oC$.
Since (9.2) does not invoke $R$, and since $\pi\gf$
is already known to satisfy (9.1) with $R=\epsilon$
(see Def.\ 4.1)
it follows that $\pi\gf\in\cQ_M(\geps)\otimes\oC$ as well.
By a second application of Lemma 9.1, this time with
$R=\epsilon$ and $\gl_0=0$, we now obtain (9.7). 
The identity (9.8) follows from (4.4) and (9.7).
The map $\gf\mapsto\pi\gf$ is continuous 
$\PACX\to\cQ\otimes\oC$ and $\Pal_M(X\col\tau)\cap\PACX$
is the preimage of $\cQ_M\otimes\oC$. Hence it is
closed.\hfq

\section{A different characterization of the Paley-Wiener space}

In [4, Def.\ 21.6], we defined the Paley-Wiener space $\PWXt$
somewhat differently from Definition 3.4, and we
conjectured in [4, Rem.~21.8], that this space
was equal to $\Fou(\CciXt)$. The purpose of this section 
is to establish equivalence of the two definitions of $\PWXt$
and to confirm the conjecture of~[4]. 

The essential difference between the definitions is that in [4]
several properties are required only 
on $\baraqd(P,0)$; the identity 
$\gf(s\gl)=\nC(s\col\gl)\gf(\gl)$ 
(cf.\ Lemma 3.10) is then part of the definition of the
Paley-Wiener space. In the following theorem we establish
a property of $\nC(s\col\gl)$ which
is crucial for comparison of the definitions. 
Let $\Pi_{\gS,\real}$ denote the set of polynomials on
$\faqdc$ which are products of functions of the form
$\gl\mapsto\langle\alpha ,\gl \rangle+c$ with $\alpha\in\gS$ and $c\in\real$.

\begin{theorem} Let $s\in W$ and let $\omega\subset\faqd$ be compact.
There exist a polynomial $q\in\Pi_{\gS,\real}$ and a number
$N\in\Natural$ such that
$\gl\mapsto(1+|\gl|)^{-N} q(\gl)\nC(s\col\gl)$ is bounded on $\omega+i\faqd$.
\end{theorem}

\Proof See [10]. 
\hfq

\begin{lem} The space $\PXt$ of Definition {\rm 4.1}
is equal to the space of
$\oC$-valued meromorphic functions on $\faqdc$
that have the properties {\rm (i)--(ii)} of Lemma {\rm 9.2}{}
together with\/{\rm :}\/
\smallskip

{\rm (iii)} For every compact
set $\omega\subset\baraqd(P,0)\cap\faqd$ and for all $n\in\Natural${\rm ,}
$$\sup_{\gl\in\omega+i\faqd} (1+|\gl|)^n \|\pi(\gl)\gf(\gl)\|<\infty.$$
Moreover{\rm ,} there exist a real $\gS$-configuration
$\Hypt$ and a map $\dt\colon\Hypt\to\Natural$ such that
\begin{equation}
\PXt\subset\Pal(\faqd,\Hypt,\dt)\otimes\oC.
\end{equation}
\end{lem}

\Proof Condition (i) in Definition 4.1 is the same as (i) in Lemma
9.2, whereas (ii) is stronger. However, it was seen in
Lemma 9.2{} that (i)$\land$(ii) implies (ii) of Definition 4.1.
The condition (iii) in Definition 4.1
is also stronger than (iii) above. 

It thus remains to be seen that (i)--(iii) above
imply (iii) of Definition 4.1, and that (10.1) holds.
We will establish both at the same time.
Let $\Hypt$ and $\dt$ be as in Lemma 9.2, and
assume that $\gf$ satisfies (i)--(iii) above; then $\gf\in
\Mer(\faqd,\Hypt,\dt)\otimes\oC$.
Let $\omega\subset\faqd$ be compact.
Using Theorem 10.1{} we see from (iii) together with (i) that 
there exists a polynomial $Q\in\Pi_{\gS,\real}$ such that
$$\sup_{\gl\in\omega+i\faqd} (1+|\gl|)^n\|Q(\gl)\gf(\gl)\|<\infty$$
for each $n\in\Natural$. Clearly we may assume that $Q$ is
divisible by $\pi_{\omega,\dt}(\gl)$ (see (2.6)). Using [2, Lemma 6.1]
and the fact that $\omega$ was arbitrary,
we can in fact remove all factors of $Q/\pi_{\omega,\dt}(\gl)$
from the estimate, so that we may assume $Q=\pi_{\omega,\dt}(\gl)$.
Hence $\gf\in\Pal(\faqd,\Hypt,\dt)\otimes\oC$. 
The statement in (iii) of Definition 4.1 follows 
by the same reasoning, when we invoke the already established statement
(ii) of that definition.
\hfq

\begin{lem} The pre-Paley-Wiener space $\Mer(X\col\tau)$
defined in {\rm [4, Def.\ 21.2],} is identical with $\cup_{M>0} \PMXt${\rm ,}
where $\PMXt$ is as defined in Definition {\rm 4.1. }
\end{lem}

\Proof Let $M>0$ and $\gf\in\PMXt$. Then properties (a) and (b)
of [4, Def.\ 21.2], are obviously fulfilled, and (c), with $R=M$, 
follows from (9.7). Hence $\gf\in\Mer(X\col\tau)$.

Conversely, let $\gf\in\Mer(X\col\tau)$, then $\gf\in\PXt$ by
Lemma 10.2. Moreover, condition
(iv) in Definition 4.1
results easily from (c) of [4], with $M=R$.
Hence $\gf\in\PMXt$.
\Endproof\vskip4pt 

In [4] the space $\PWXt$ is defined as the space of functions
$\gf\in\Mer(X\col\tau)$ that satisfy certain relations. These relations
will now be interpreted in terms of Laurent functionals by means of the
following lemma.

\begin{lem} Let $u_1,\dots,u_k\in S(\faqd)${\rm ,}
$\psi_1,\dots,\psi_k\in\oC${\rm ,} and $\gl_1,\dots,\gl_k\in\baraqd(P,0)$. 
Then there exists a Laurent functional $\cL\in\LauoC${\rm ,} such that
\begin{equation}
\cL\gf=\sum_{i=1}^k 
u_i[\pi(\gl)\hinp{\gf(\gl)}{\psi_i}]_{\gl=\gl_i}
\end{equation}
for all $\gf\in\Mer(X\col\tau)$. 
Conversely{\rm ,} given $\cL\in\LauoC$
there exist $k${\rm ,} $u_i${\rm ,} $\psi_i$ and $\gl_i$ as above such that\/
{\rm (10.2)} holds for all $\gf\in\Mer(X\col\tau)$. 
\end{lem}

\Proof To prove the existence of $\cL$
we may assume that $k=1$. Let $d=d_{X,\tau}$ and let
$\pi_1=\pi_{\{\gl_1\},d}$ be determined by (2.6). Then $\pi_1$
divides $\pi$; let $p$ denote their quotient.
It follows from [7, Lemma~10.5], 
that there exists $\cL_1\in\LauoC$ such that 
$$\cL_1\gf=u_1[\pi_1(\gl)\hinp{\gf(\gl)}{\psi_1}]_{\gl=\gl_1}$$
for all $\gf$ such that $\pi_1\gf$ is holomorphic
near $\gl_1$. By Lemma 7.3  the map $\cL\maps\gf\mapsto\cL_1(p\gf)$ 
belongs to $\LauoC$. It clearly satisfies (10.2).

Conversely, let $\cL\in\LauoC$ be given. We may assume that the
support of $\cL$ consists of a single point in $\faqcd$.
This point equals $s\gl_0$ for suitable
$\gl_0\in\baraqd(P,0)$ and $s\in W$. 
Let $\pi_0$, $\pi_s$ and $\pi^\sim$ be as in the
proof of Lemma 9.2.
The restriction of $\cL$ to $\Mer(\faqd,\Hypt,\dt)\otimes\oC$
is a finite sum of terms of the form
\begin{equation}
\gf\mapsto u[\pi^\sim(\gl)\hinp{\gf(\gl)}{\psi}]_{\gl=s\gl_0},
\end{equation}
where $\psi\in\oC$ and $u\in S(\faqd)$. 
For $\gf\in\Mer(X\col\tau)$ we use the Weyl conjugation property 
and rewrite (10.3) in the form
$$\gf\mapsto u[\pi^\sim(s\gl)
\hinp{\nC(s\col\gl)\gf(\gl)}{\psi}]_{\gl=\gl_0},$$
in which the element $u$ has been replaced by its $s$-conjugate.
Since the polynomial $\pi_0\pi_s$ divides $\pi^\sim(s\gl)$, 
and since $\pi_s(\gl)\nC(s\col\gl)$ is holomorphic at $\gl_0$
it follows from the Leibniz rule that this expression
can be further rewritten as a finite sum of terms of the form
\begin{equation}
\gf\mapsto u[\pi_0(\gl)
\hinp{\gf(\gl)}{\psi}]_{\gl=\gl_0}
\end{equation}
where $\psi\in\oC$ and $u\in S(\faqd)$. 
Finally, since $\pi_0$ divides $\pi$, 
the following lemma shows that there exists $u'\in S(\faqd)$
such that (10.4) takes the form
$\gf\mapsto u'[\pi(\gl)
\hinp{\gf(\gl)}{\psi}]_{\gl=\gl_0},$
which is as desired in (10.2).
\Endproof\vskip4pt 

Let $\Pi_\real$ denote the set of polynomials on
$\faqdc$ which are products of functions of the form
$\gl\mapsto  \inp{\xi}{\gl} +c$ with 
$\xi \in \faqd\setminus\{0\}$ and $c\in\real$.

\begin{lem} Let $p\in\Pi_\real$. There exists {\rm } for each $u\in S(\faqd)${\rm ,}
an element $u'\in S(\faqd)$ such that $u'(p\gf)(0)=u\gf(0)$ for all
germs $\gf$ at $0$ of holomorphic functions on $\faqcd$.\end{lem}

\Proof 
We may assume that the degree
of $p$ is one. Then $p(\gl)=\inp{\xi}{\gl}+p(0)$ for some nonzero
$\xi\in\faqd$. The case that $p(0) = 0$ is covered by 
[5, Lemma 1.7 (i)]. Thus, we may assume that $p(0) =1.$
Let $\xi'=\xi/\inp{\xi}{\xi}$. Then
$\xi' p=1$, when $\xi'$ is considered as a constant 
coefficient differential operator acting on the function $p$.
By linearity we may assume that $u$ is of the form
$u=u''\xi'{}^k$ with $k\in\Natural$ and $u''\in S(\xi'{}^\perp)$.
Let
$u'=u''\sum_{i=0}^k (-1)^{k-i}\frac{k!}{i!} \xi'{}^i.$
A simple calculation with the Leibniz rule shows that
$u'(p\gf)(0)=u\gf(0)$, as desired.
\hfq

\begin{cor}  The Paley-Wiener spaces\/ $\PWXt$ in Definition {\rm 3.4} 
and in {\rm [4, {\it Def.}\ 21.6],} are identical{\rm ,} and both are equal to
$\Fou(\CciXt))$.
\end{cor}

\Proof In view of (9.7), 
it is immediate from Lemmas 10.3{} and 10.4{}
that the space $\PWXt$ of [4] is identical to the space
denoted $\PWC$ in Remark 3.9. According to that remark,
it follows from Theorem 3.6{} that this space is equal to
$\PWXt$ as well as to $\Fou(\CciXt)$.
\hfq

\end{document}